\documentclass[12pt]{article}
\usepackage{color}
\usepackage{amsmath,amsfonts,amssymb,graphics,amsthm}
\usepackage{hyperref}
\usepackage {graphicx}
\usepackage {verbatim}
\usepackage[margin=1in]{geometry}
\usepackage[normalem]{ulem}

\newtheorem{theorem}{Theorem}[section]
\newtheorem{corollary}[theorem]{Corollary}
\newtheorem{lemma}[theorem]{Lemma}
\newtheorem{proposition}[theorem]{Proposition}

\newtheorem{remark}[theorem]{Remark}

\newcommand{\old}[1]{{}}

\def\MR#1{\href{http://www.ams.org/mathscinet-getitem?mr=#1}{MR#1}}

\newcommand{\aryb}{\begin{eqnarray*}}
\newcommand{\arye}{\end{eqnarray*}}
\def\alb#1\ale{\begin{align*}#1\end{align*}}
\newcommand{\eqb}{\begin{equation}}
\newcommand{\eqe}{\end{equation}}
\newcommand{\eqbn}{\begin{equation*}}
\newcommand{\eqen}{\end{equation*}}

\newcommand{\C}{\mathbb{C}}
\newcommand{\D}{\mathbb{D}}
\newcommand{\E}{\mathbb{E}}
\newcommand{\N}{\mathbb{N}}
\newcommand{\Q}{\mathbb{Q}}
\newcommand{\Z}{\mathbb{Z}}
\newcommand{\R}{\mathbb{R}}
\renewcommand{\P}{\mathbb{P}}

\newcommand{\eps}{\varepsilon}
\newcommand{\wh}{\widehat}
\newcommand{\wt}{\widetilde}
\newcommand{\1}{\mathbf{1}}

\DeclareMathOperator{\SLE}{SLE}

\def\cT{\mathcal{T}}
\def\cS{\mathcal{S}}

\def\cO{\mathcal{O}}

\def\cL{\mathcal{L}}

\def\cC{\mathcal{C}}

\newcommand{\BB}{\mathbb}
\newcommand{\mcl}{\mathcal}
\newcommand{\ol}{\overline}
\newcommand{\ul}{\underline}

\newcommand{\eqd}{\overset{d}{=}}

\newcommand{\arxiv}[1]{\href{http://arxiv.org/abs/#1}{#1}}

\newcommand{\op}{\operatorname}

\newcommand{\frk}{\mathfrak}
\newcommand{\eqD}{\overset{d}{=}}
\newcommand{\ep}{\epsilon}
\newcommand{\rta}{\rightarrow}

\begin {document}
\author{Nina Holden\footnote{Massachusetts Institute of Technology, Cambridge, MA, ninah@math.mit.edu}
	\qquad\qquad Xin Sun\footnote{Columbia University, New York, NY, xinsun@math.columbia.edu} }
\title{SLE as a mating of trees in Euclidean geometry}
\date{}
\maketitle

\begin{abstract}
The mating of trees approach to Schramm-Loewner evolution (SLE) in the random geometry of Liouville quantum gravity (LQG) has been recently developed by Duplantier-Miller-Sheffield (2014). In this paper we consider the mating of trees approach to SLE in Euclidean geometry. Let $\eta$ be  a whole-plane space-filling SLE with parameter $\kappa>4$, parameterized by Lebesgue measure. The main observable in the mating of trees approach is the \emph{contour function},
a two-dimensional continuous process describing the evolution of the Minkowski content of the left and right frontier of $\eta$.
We prove regularity properties of the contour function and show that (as in the LQG case) it encodes all the information about the curve $\eta$. We also prove that the uniform spanning tree on $\Z^2$ converges to $\SLE_8$  in the natural topology associated with the mating of trees approach.
\end{abstract}
\tableofcontents

\section{Introduction}
The Schramm-Loewner evolution (SLE) is a one-parameter family of random fractal curves introduced by Oded Schramm as a candidate for scaling limits of interfaces in two-dimensional statistical physics models \cite{schramm0}. Since it was introduced, SLE has proved to be the limit of several lattice models, see e.g.\ \cite{lsw-lerw-ust,smirnov-cardy,ss-dgff,chelkak-smirnov-ising,cdhks14,ks-ising,lawler-viklund-16}.

\begin{figure}[ht!]
	\begin{center}
		\includegraphics[scale=2]{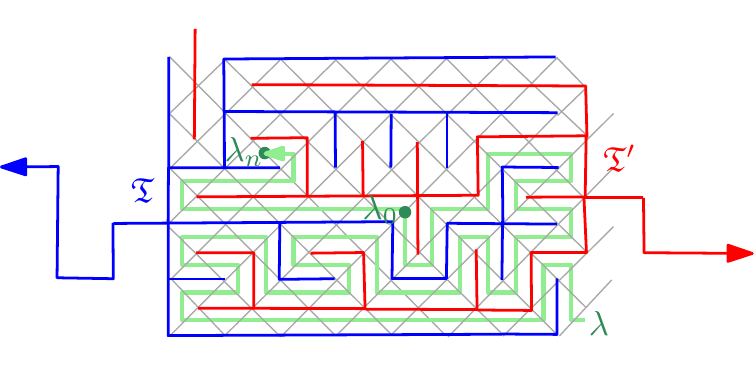} 
		\caption{A spanning tree $\frk T$ on $\Z^2$ (blue), its dual tree $\frk T'$ (red), and the Peano curve $\lambda$ (green). The Peano curve traces the interface of $\frk T$ and $\frk T'$ at unit speed, meaning that it takes one unit of time to traverse each gray triangle. The pair of functions $(\frk L,\frk R)$ encodes the height in the pair of trees $(\frk T,\frk T')$, such that for each $n\in\Z$, $\frk L_n$ (resp.\ $\frk R_n$) denotes the height in $\frk T$ (resp.\ $\frk T'$) at position $\lambda_n$, relative to the height in the tree at position $\lambda_0$. 
			The blue (resp.\ red) arrow points to the root of $\frk T$ (resp.\ $\frk T'$) at $\infty$. 
		} \label{fig-ust}
	\end{center}
\end{figure}

Given a uniform spanning tree (UST) $\frk T$ on $\Z^2$, there is a.s.\ a uniquely determined spanning tree $\frk T'$ in the dual graph, which is defined such that $\frk T$ and $\frk T'$ never cross each other, see Figure \ref{fig-ust}. The Peano curve $\lambda$ is the interface between $\frk T$ and $\frk T'$.   It was proved in \cite{lsw-lerw-ust} that in a chordal setting the Peano curve $\lambda$ of a uniform spanning tree converges in law in the scaling limit to an SLE$_8$ $\eta$ in the space of curves equipped with the $L^\infty$ norm, viewed modulo reparametrization of time.

Throughout this paper we define $\lambda$ as follows (see Figure \ref{fig-ust}). We let $\lambda$ be a function from $\R$ to $\C$ with $\lambda_0=(\tfrac14,\tfrac14)$, and such that for all $t\in \Z$, $\lambda|_{[t,t+1]}$ is a straight line segment of length $\tfrac12$ in up, down, left or right direction. Moreover, $\lambda$ is the interface between $\frk T$ and $\frk T'$ so that $\frk T$ is on the left side of $\lambda$. For each $n\in\Z$, the point $\lambda_n$ is contained in the line segment between points $(k_n,m_n)\in\Z^2$ and $(k'_n,m'_n)\in(\Z+\frac 12)^2$ satisfying $|k_n-k'_n|=|m_n-m'_n|=\frac 12$. Let $(\wh k_n,\wh m_n)\in \Z^2$ be the first point on the path from $(k_n,m_n)$ to $\infty$ in $\frk T$ which is also on the path from $(0,0)$ to $\infty$. Let $\frk L_n$ be the $\frk T$-graph distance from  $(k_n,m_n)$ to $(\wh k_n,\wh m_n)$, minus the $\frk T$-graph distance from $(0,0)$ to $(\wh k_n,\wh m_n)$. We define $\frk R_n$ similarly by considering $\frk T'$ instead of $\frk T$.
 We say that $\frk Z=(\frk L,\frk R)$ \emph{encodes} the trees $\frk T$ and $\frk T'$, since, as we will explain later, $\frk T$ and $\frk T'$ are measurable with respect to $\frk Z$ up to rotation by $\frac{\pi}{2}$ about the origin.

We can define the corresponding contour functions $Z=(L_t,R_t)_{t\in\R}$ for the continuum scaling limit $\eta$, which is an SLE$_8$ in $\C$ from $\infty$ to $\infty$. Let $\eta$ be parametrized by Lebesgue measure, i.e., if $\mcl L$ denotes Lebesgue measure then $\mcl L(\eta([s,t]))=t-s$ for any $s<t$, and let $\eta(0)=0$. Given an enumeration $(z_n )_{n\in\N}$ of $\Q^2$, for each $n\in\N$ let $\eta^L_{z_n}$ (resp.\ $\eta^R_{z_n}$) be the curve describing the left (resp.\ right) frontier of $\eta$ stopped upon hitting $z_n$. By SLE duality these curves have the law of whole-plane SLE$_2$. The set of curves $\{\eta^L_{z_n}\,:\, n\in\N \}$ defines a space-filling tree $\mcl T$, where each curve $\eta^L_{z_n}$, $n\in\N$, is a branch of $\mcl T$ from the leaf $z_n$ to the root of $\mcl T$ at $\infty$. Similarly, the set of curves $\{\eta^R_{z_n}\,:\, n\in\N \}$ defines a dual space-filling tree $\mcl T'$, and it is immediate from the construction that the branches of $\mcl T$ and $\mcl T'$ never cross each other. As we will explain in more detail later, by properties of the natural parametrization of SLE \cite{lawler-shef-nat,lawler-minkowski,LV-radial}, the natural length measure along the branches of $\mcl T$ and $\mcl T'$ is the $5/4$-dimensional Minkowski content of the curves $\eta_{z_n}^L$ and $\eta_{z_n}^R$.  Let $L_t$ (resp.\ $R_t$) denote the height in $\mcl T$ (resp.\ $\mcl T'$) at time $t\in\R$, relative to the height in $\mcl T$ (resp.\ $\mcl T'$) at time 0, when we use the Minkowski content to measure the length of the branches. 

Our first result is that $Z$ is well-defined, and is the scaling limit of $\frk Z$. Consider an instance of the UST on $\Z^2$  and the associated Peano curve $\lambda$. 
For all $\delta\in(0,1]$,  let $\eta^\delta(t)=: \delta\lambda_{\delta^{-2}t}$. For $t\in\delta^2\Z$ define $L^\delta_t:=\check c\delta^{5/4}\frk L_{\delta^{-2}t}$, where $\check c>0$ is a universal constant (which is the same as the one appearing in Theorem \ref{thm:lawler-viklund-wholeplane}), and for $t\not\in\delta^2\Z$ define $L^\delta_t$ by linear interpolation. The function $R^\delta$ is defined similarly.  We view $\eta$ and $\eta^\delta$ as elements in the set of parametrized curves on $\C$ equipped with the topology of uniform convergence on compact sets. The contour functions $Z$ and $Z^\delta$ are elements in the space of two-dimensional continuous functions equipped with the topology of uniform convergence on compact sets.

\begin{theorem}\label{thm:convergence}
	For $\delta\in(0,1]$, consider a UST  on $\delta\Z^2$ and an instance of a whole-plane space-filling SLE$_8$ $\eta$ in $\C$. With the notation introduced above, $Z=(L,R)$ is well-defined as a continuous function, and the pair
	$(\eta^\delta,Z^\delta)$ converges in law to $(\eta,Z)$ as $\delta\rta 0$.
\end{theorem}

\begin{remark}
	Theorem \ref{thm:convergence} implies that the UST and its dual also converge in the space whose elements are measured, rooted real trees continuously embedded into $\C$ (see \cite[Section 3]{barlow-sub} for the precise definition of this topology). Tightness of the UST in this topology was proved in \cite{barlow-sub}. The convergence result follows from the above theorem, since the functions $L^\delta$ and $R^\delta$ are rescaled version of the UST and dual tree contour functions (up a time change of $o_\delta(1)$), and since convergence of contour functions implies convergence in the Gromov-Hausdorff-Prokhorov topology (see e.g.\ \cite[Proposition 3.3]{ghp13}).
\end{remark}

We may proceed similarly as above to define contour functions $Z=(L_t,R_t)_{t\in\R}$ for SLE$_\kappa$ for other values of $\kappa$.
Let $\kappa>4$, and let $\eta$ be a whole-plane space-filling SLE$_\kappa$ $\eta$ from $\infty$ to $\infty$ as defined in Section \ref{sec:ig}. Similarly as above, we let $\eta$ be parametrized by Lebesgue measure and satisfy $\eta(0)=0$, and for any $z\in\C$ let $\eta_z^L$ (resp.\ $\eta_z^R$) denote the left (resp.\ right) frontier of $\eta$ when the curve first hits $z$. Given any $t\in\R$ let $L_t$ (resp.\ $R_t$) denote the length of $\eta_{\eta(t)}^L$ (resp.\ $\eta_{\eta(t)}^R$) relative to the length of $\eta_{0}^L$ (resp.\ $\eta_{0}^R$). Lengths are measured by considering the natural parametrization of the curves, which is given by $(1+2/\kappa)$-dimensional Minkowski content.

\begin{figure}[ht!]
	\begin{center}
		\includegraphics[scale=1]{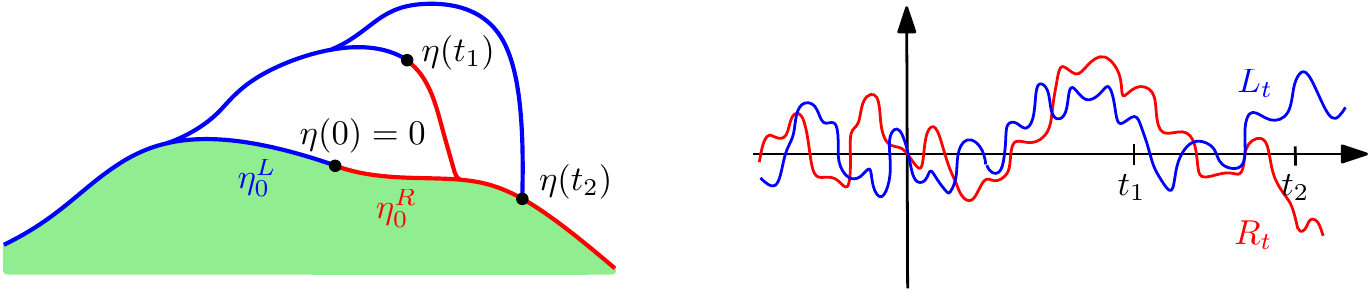} 
		\caption{For $\kappa>4$ and an SLE$_\kappa$ $\eta$, the function $Z=(L,R)$ describes the evolution of the left and the right, respectively, boundary length of $\eta$. The boundary length is measured in $(1+2/\kappa)$-dimensional Minkowski content. The time $t_2>0$ on the figure is a time at which $R$ reaches a running infimum relative to time 0. We remark that $\eta((-\infty,0])$, which is shown in green, has a different topology than on the figure for $\kappa\in(4,8)$. } \label{fig-Z}
	\end{center}
\end{figure}

In part $(iii)$ of the theorem below we let $\cC(\R,\R^2)$ denote the space of equivalence classes of continuous processes $W=(W_t)_{t\in\R}$ with values in $\R^2$, such that $W^1$ and $W^2$ are equivalent if there exists an increasing bijection $s:\R\to\R$ such that for all $t\in\R$, we have $W^2_{t}=W^1_{s(t)}$.
\begin{theorem}
	Let $\kappa>4$, and let $\eta$ and $Z$ be as above.
	\begin{itemize}
		\item[(i)] The process $Z$  is a.s.\ well-defined as an $\alpha$-H\"{o}lder continuous process for any $\alpha<1/2+1/\kappa$, and the following probability decays faster than any power of $M$ for fixed $\alpha$
		\eqb
		\P\left[ \sup_{s,t\in[0,1]} \frac{|Z_t-Z_s|}{|s-t|^\alpha}>M \right].
		\eqe 
		\item[(ii)] For any $a>0$, $(Z_{t})_{t\in\R}\overset{d}{=}(a^{1/2+1/\kappa}Z_{a^{-1}t})_{t\in\R}$. The process $Z$ has stationary increments, and the tail $\sigma$-algebra of $Z$ is trivial. Furthermore,
		\eqb
		\limsup_{t\rta\pm\infty} V_t = \infty,\quad
		\liminf_{t\rta\pm\infty} V_t = -\infty,\quad\text{for\,\,}V=L,R.
		\label{eq:intro1}
		\eqe
		\item[(iii)] Assume $\kappa=8$. The process $Z$ defines an object $Z'$ in the space $\cC(\R,\R^2)$. It holds that $Z'$ determines $Z$, i.e., $Z$ is measurable with respect to the $\sigma$-algebra generated by $Z'$.
	\end{itemize}
	\label{thm:Z}
\end{theorem}
Theorem \ref{thm:Z} will be proved in Section \ref{sec:contour}. 
In the proof we use the mating of trees theorem in the Liouville quantum gravity setting (see below) to deduce the desired properties of the contour functions in the Euclidean setting. The reason we only prove $(iii)$ for the case $\kappa=8$, is that we need a lower bound for the Minkowski content of the frontier, which we only know for $\kappa=8$, although we expect it to hold also for other $\kappa$.
A more substantial part of the paper  is devoted to proving the following theorem, asserting that the mating of trees in Euclidean geometry encodes all the information of the space-filling SLE.  
Hence the mating of trees provides an alternative way of encoding conformal invariant systems other than interfaces which have SLE as their scaling limits.
The proof is given in Section \ref{sec:graph}, using results from Sections \ref{sec:ig}, \ref{sec:contour} and \ref{sec:expectation}. The proof crucially relies on the assumption that the shortest path between two points is the straight line, a defining property of Euclidean geometry. (See Proposition~\ref{prop:philinear}.)
Another technical ingredient is a regularity estimate for space-filling SLE (see Proposition~\ref{prop:powerlaw})  proved via imaginary geometry, which is of independent interest.  
\begin{theorem}
	Let $\eta$ be a whole-plane space-filling SLE$_\kappa$ for $\kappa>4$ in $\C$, and define $Z$ as in Theorem \ref{thm:Z}. Then $\eta$ is measurable with respect to the $\sigma$-algebra generated by $Z$, modulo rotations of $\eta$ about the origin.
	\label{thm:meas}
\end{theorem}
The analogous result to this theorem in  the context of Liouville quantum gravity (LQG)  was proved in \cite{wedges}, see further details after Corollary~\ref{cor:mate}. Our proof is  different in nature  as it relies on the Euclidean geometry. The discrete analogue of Theorem \ref{thm:meas} for $\kappa=8$ says that a spanning tree $\frk T$ on $\Z^2$ is measurable with respect to the pair of contour functions $(\frk L,\frk R)$ of $\frk T$ and its dual $\frk T'$ up to a $\tfrac{\pi}{2}$-rotation. This discrete result follows from e.g.\ a bijection of Mullin \cite{mullin-maps} (see also \cite{bernardi07,burger}) in the context of planar maps.

The result that $(L^\delta,R^\delta)$ converges in law to $(L,R)$ means that the UST on $\Z^2$ converges to SLE$_8$ in a Euclidean analogue of the mating of trees topology, which was used in \cite{burger,burger1,burger2,burger3,kmsw-bipolar,gkmw-burger,ghs-bipolar,Schnyder} to prove convergence of decorated random planar maps to SLE-decorated LQG. There  tree-decorated discrete models are said to converge to SLE-decorated LQG in the mating of trees sense if the contour functions of the trees converge to a pair of correlated Brownian motions encoding a pair of continuum random trees.  (See more discussion below Corollary~\ref{cor:mate}.)

The natural parametrization of SLE$_\kappa$ is a parametrization which is conjectured (or proved, for $\kappa=2$) to capture the natural parametrization of the associated discrete models, i.e., one unit of time corresponds to traversing one edge/vertex/face of the discrete model. It is therefore natural to conjecture that $Z$ for other values of $\kappa$ is the scaling limit of the contour functions of other discrete tree-decorated models having SLE$_\kappa$ as a scaling limit. We remark that such convergence results would follow by proceeding as in Section \ref{sec:tight}, once analogues of Theorem \ref{thm:lawler-viklund-wholeplane}, Proposition \ref{prop:lsw} and Lemma \ref{lem:metric-compare} (by other authors) were established. Certain discrete models which are conjectured to converge to SLE$_\kappa$ with $\kappa>8$, for example the 6-vertex model \cite{6V} and the 20-vertex model \cite{Schnyder}  are naturally decorated with \emph{multiple} pairs of trees, and one may then hope to establish joint convergence of these trees by proving joint convergence of the corresponding pairs of contour functions, similarly to the results established in \cite{ghs-bipolar} for random planar maps.

\begin{figure}[ht!]
	\begin{center}
		\includegraphics[scale=0.8]{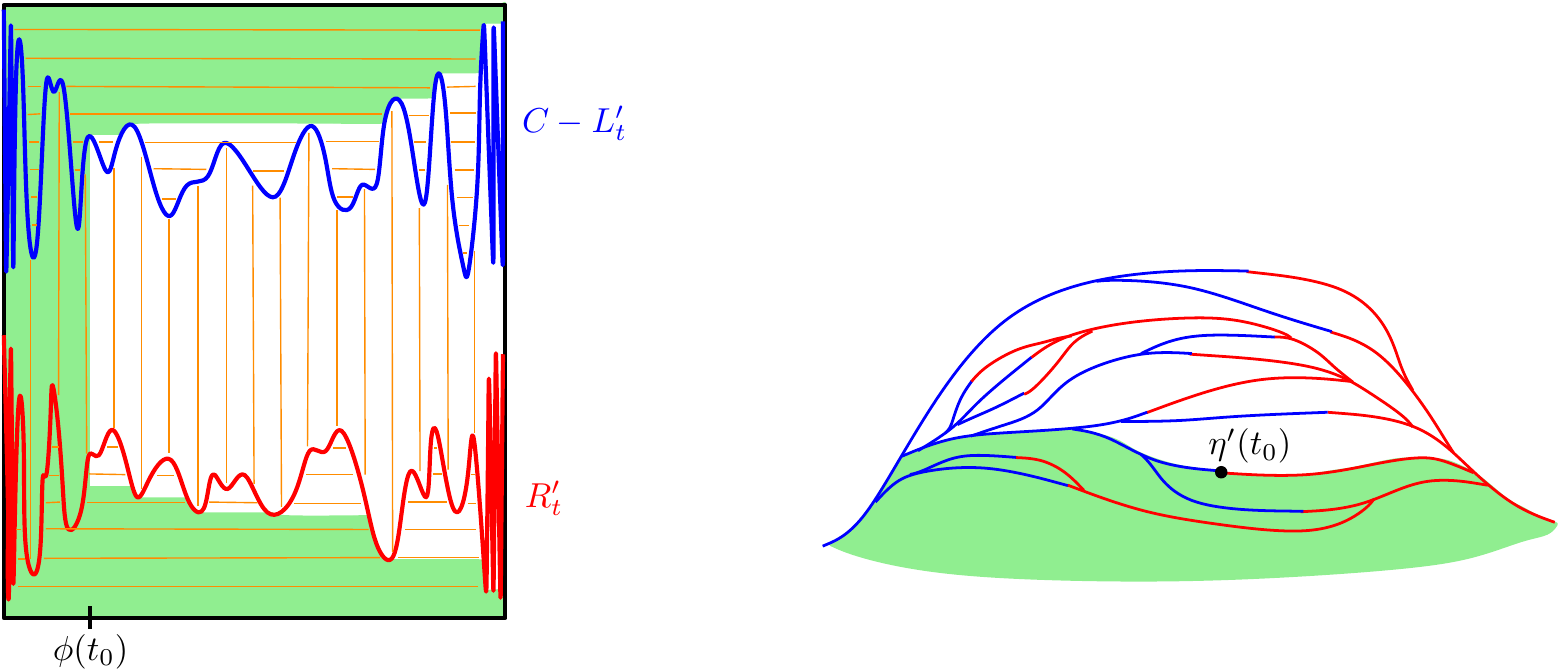} 
		\caption{The figure illustrates how we obtain a topological sphere decorated with a space-filling path from a pair of functions $L,R$ satisfying \eqref{eq:intro1}. Each function $L,R$ encodes an infinite tree (shown in blue and red, respectively, on the right figure), and the idea of the construction is to glue together these two trees. Letting $\phi:\R\to(0,1)$ be a strictly increasing bijective map we let $(L'_t)_{t\in(0,1)}$ and $(R'_t)_{t\in(0,1)}$ be $(0,1)$-valued processes defined by $L'_t:=\phi(L_{\phi^{-1}(t)})$ and $R'_t:=\phi(R_{\phi^{-1}(t)})$. For some constant $C>0$ we draw $R'$ and $C-L'$ in a rectangle as on the left figure, where $C$ is chosen sufficiently large such that the two curves don't intersect. We define an equivalence relation on the rectangle by identifying (i) all points on the boundary of the rectangle, (ii) all points that lie on the same line segment below $R'$ (resp.\ above $L'$), and (iii) all points that lie on the same vertical line between $R'$ and $C-L'$. We will argue (inspired by arguments in \cite{wedges}) that the set of equivalence classes just defined gives a topological sphere. The sphere is decorated with the space-filling path which maps $t_0\in\R$ to the equivalence class of the point $(\phi(t_0),R'_{\phi(t_0)})$. This figure first appeared in \cite{ghs-bipolar}.	
		} \label{fig-mating}
	\end{center}
\end{figure}

Since $L$ and $R$ are continuous functions satisfying \eqref{eq:intro1}, the functions $L$ and $R$ are the contour functions of a pair of infinite-volume real trees \cite{legall-trees-survey}. Inspired by \cite{wedges}, we deduce from Theorem \ref{thm:meas} that we may ``glue'' together the two trees to obtain a topological sphere decorated with a space-filling path, which can then be embedded canonically into the complex plane. See Section \ref{sec:graph} for a proof of the following corollary, and see Figure \ref{fig-mating} for an illustration.
\begin{corollary}\label{cor:mate}
	For $\kappa>4$ and $Z$ with the same marginal law as in Theorem \ref{thm:meas}, we obtain a topological sphere with a space-filling path when gluing together the associated pair of trees as explained in Figure \ref{fig-mating}. This path-decorated sphere has a canonical embedding into the complex plane, where the law of the curve is that of a space-filling SLE$_\kappa$.
	\label{cor1}
\end{corollary}

Finally we will describe an analogue of Theorem \ref{thm:meas} and its corollary in the context of LQG  \cite{wedges}. In this setting $Z$ has the law of a two-dimensional correlated Brownian motion. The curve $\eta$ still has the law of a space-filling SLE$_\kappa$, $\kappa>4$, but it lives on top of a $\gamma$-LQG surface ($\gamma=4/\sqrt{\kappa}$) which determines the parametrization of $\eta$ and induces a measure on the frontier of $\eta$.

Recall that for any $\gamma\in(0,2)$ and a domain $D\subset\C$, $\gamma$-Liouville quantum gravity \cite{shef-kpz,rhodes-vargas-review} is a random surface which may be written heuristically as $e^{\gamma h}\,dz$, where $h$ is an instance of a Gaussian free field (GFF) \cite{shef-gff} or a related form of distribution in $D$ and $dz$ denotes Lebesgue measure in $D$. The term $e^{\gamma h}$ does not make literal sense since $h$ is a distribution and not a function, but as explained in the above references it has been made sense of as a random area measure in $D$. The GFF also induces a random length measure along certain curves in $D$. 

For any $\kappa>4$ and $\gamma:=4/\sqrt{\kappa}$ the authors of \cite{wedges} considered a pair of Brownian motions $Z=(L_t,R_t)_{t\in\R}$ with correlation $-\cos(4\pi/\kappa)$ satisfying $(L_0,R_0)=(0,0)$ (see \cite{cov8} for the correlation when $\kappa>8$). By ``gluing'' together the corresponding infinite volume continuum random trees \cite{aldous-crt1,aldous-crt2,aldous-crt3} as in Figure \ref{fig-mating}, they obtained a topological sphere with a space-filling path and an area measure, called a \emph{peanosphere}. They then proved an analogue of Corollary \ref{cor:mate} above, namely that the peanosphere has a canonical embedding into $\C$, where the space-filling path has the law of an SLE$_\kappa$ $\eta$, and the area measure has the law corresponding to an independent instance of the $\gamma$-LQG surface known as the $\gamma$-quantum cone \cite[Section 4.2]{wedges}.

Alternatively, their result can be stated as in the following theorem. Consider a space-filling SLE$_\kappa$ $\eta$ which lives on an independent $\gamma$-quantum cone with area measure $\mu$. Parametrize $\eta$ by $\gamma$-LQG area measure, i.e., $\mu(\eta([s,t]))=t-s$ for any $s<t$, and let $\eta(0)=0$. The $\gamma$-LQG surface defines a length measure along the frontier of $\eta((-\infty,t])$ at any fixed time $t\in\R$. Let $L_t$ (resp.\ $R_t$) denote the length of the left (resp.\ right) frontier of $\eta((-\infty,t])$ relative to the length at time 0. Set $Z=(L,R)$. The following is \cite[Theorem 1.13]{wedges}.
\begin{theorem}[\cite{wedges}]
	In the setting above, $(\eta,h)$ is measurable with respect to the $\sigma$-algebra generated by $Z$.
	\label{thm:wedges}
\end{theorem}

In Section \ref{sec:ig} we review imaginary geometry and the construction of space-filling SLE, and we prove some basic lemmas which are needed in the remainder of the paper. Theorem \ref{thm:convergence} and Theorem \ref{thm:Z} are proved in Section  \ref{sec:tight} and Section \ref{sec:contour}, respectively. In Section \ref{sec:graph} we prove Theorem \ref{thm:meas} and Corollary \ref{cor1}, modulo two technical results which are proved in Section \ref{sec:expectation}.

\subsection{Notation}
We write $a\preceq b$ (resp.\ $a\succeq b$) if there is a constant $C$ independent of the parameters of interest such that $a\leq Cb$ (resp.\ $a\geq Cb$). We write $a\asymp b$ if $a\preceq b$ and $a\succeq b$. We say that $f(n)$ has \emph{superpolynomial decay} if  $f(n)\preceq n^{-p}$ for any $p$ as $n\to \infty$.

For any $z\in\C$ and $r>0$ we let $B_r(z):=\{w\in\C \,:\, |z-w|<r\}$ be the Euclidean ball of radius $r$ centered at $z$. We let $\D=B_1(0)$ be the unit disk centered at the origin. 

For any $D\subset\C$ we let $\frk m(D)$ denote the $d$-dimensional Minkowski content of $D$, where the dimension $d$ is implicitly understood to be given by $d=1+2/\kappa$ when we work with SLE$_\kappa$ or SLE$_{16/\kappa}$ for $\kappa>4$. Throughout the paper we will use $\ul\kappa$ (rather than $\kappa$) when we consider SLE parameters smaller than 4, and we will let $\ul\eta$ denote an associated SLE$_{\ul\kappa}$. We let $\mcl L(D)$ denote the Lebesgue measure of $D$, and we let $\op{diam}(D)$ denote the diameter of $D$.

We abuse notation in the following way throughout the paper for an arbitrary random variable $X$. When we say ``measurable with respect to $X$'' we mean ``measurable with respect to the $\sigma$-algebra generated by $X$''.

\subsubsection*{Acknowledgements}
N.H.\ was supported by a fellowship from the Norwegian Research Council. X.S. was partly supported by NSF grant DMS-1209044 and by Simons Society of Fellows. Part of this work was carried out during the Random Geometry semester at the Isaac Newton Institute, Cambridge University, and the authors would like to thank the institute and the organizers of the program for their hospitality.  The authors would also like to thank Martin Barlow, Stephane Benoist,  Ewain Gwynne, Greg Lawler, Jean-Francois Le Gall, Scott Sheffield, Wendelin Werner, David Wilson, and Dapeng Zhan for helpful discussions. They also thank the anonymous referee for his/her careful reading and for numerous helpful comments.

\section{Imaginary geometry and space-filling SLE}\label{sec:ig}
In this section we give a brief review of imaginary geometry \cite{IG1,IG2,IG3,IG4}, the construction of space-filling SLE\old{$_\kappa$ for $\kappa>4$}, and prove a few basic lemmas which will be needed later. Throughout this section and in the rest of the paper, we set
$\kappa>4$ and define
\eqb
\ul\kappa=\frac{16}{\kappa}\in (0,4), \qquad\chi = \frac{2}{\sqrt {\ul\kappa}} - \frac{\sqrt{\ul\kappa}}{2}, \qquad
\lambda = \frac{\pi}{\sqrt{\ul\kappa}},\qquad
\lambda' = \frac{\pi}{\sqrt{\kappa}}=  \frac{\pi\sqrt{\ul\kappa}}{4}.
\label{eq92}
\eqe
Let $D\subseteq \C$ be a domain and $h$ be an instance of the Gaussian free field \cite{shef-gff,IG1} in $D$. We view $h$ as a field modulo a global additive multiple of $2\pi\chi$, see \cite{IG1,IG4}.
For any given $z\in D$ and angle $\theta\in[0,2\pi)$, imaginary geometry provides a way to define the flow line $\eta^\theta_z$ for $h$ of angle $\theta$ started at $z$. The flow line may be interpreted as a solution to the following formal ODE with initial condition $\eta_z^\theta(0)=z$
\eqbn
\frac{d}{dt} \eta^\theta_z (t) = e^{ i\big( h(\eta^\theta_z(t))/\chi + \theta \big) },\qquad t>0.
\eqen
This ODE does not make literal sense, since $h$ is a distribution and not a function, but has been made sense of in \cite{IG4} (see also the earlier works \cite{dubedat-coupling,IG1} for the case $z\in\partial D$). Given an instance of $h$, imaginary geometry defines a collection of coupled flow lines $\eta_z^\theta$, simultaneously for any $z\in\Q^2$ and $\theta$ a rational multiple of $\pi$. If $\theta=\pi/2$ (resp.\ $\theta=-\pi/2$) we say that the flow line is west-going (resp.\ east-going), and we denote it by $\eta_z^L$ (resp.\ $\eta_z^R$). 

For two domains $D$ and $\wt D$, a conformal transformation $\psi:\wt D\to D$, and a field $h$ (resp.\ $\wt h$) in $D$ (resp.\ $\wt D$), we say that the pairs $(D,h)$ and $(\wt D,\wt h)$ are equivalent iff
\eqb
(\wt D,\wt h) = (\psi^{-1}(D),h\circ\psi-\chi\op{arg}\psi').
\label{eq90}
\eqe
If $\eta_{z}^\theta$ is a flow line for the field $h$ in $D$, then the image $\wt\eta_{\wt z}^{\theta}$ of $\eta_{z}^\theta$ under $\psi^{-1}$ is a flow line for $\wt h$ in $\wt D$.

In the case when $D=\C$, the marginal law of $\eta_z^\theta$ for any $\theta\in[0,2\pi)$, is that of a whole-plane SLE$_{\ul\kappa}(2-\ul\kappa)$ from $z$ to $\infty$ \cite[Theorem 1.1]{IG4}. If $\eta$ is a curve in $D$, we say that a field $h$ has \emph{flow line boundary data} if the boundary data on the left (resp.\ right) side of $\eta$ are given by $-\lambda'$ (resp.\ $\lambda'$), plus $\chi$ times the winding of the curve in counterclockwise direction. See \cite[Figure 1.9]{IG1} and \cite[Figure 1.9]{IG4}. For any stopping time $\tau$ for the flow line $\eta_z^\theta$ in $\C$, the conditional law of $h$ given $\eta_z^\theta([0,\tau])$ is that of a GFF in $\C\setminus \eta_z^\theta([0,\tau])$ with flow line boundary data. If $D\neq \C$ is a domain with harmonically non-trivial boundary and $z\in D$, then the marginal law of $\eta_z^\theta$ depends on the boundary data of $h$. For any stopping time $\tau$ for $\eta_z^\theta$, the conditional law of $h$ given $\eta_z^\theta([0,\tau])$ is that of a GFF in $D\setminus \eta_z^\theta([0,\tau])$ with flow line boundary data on $\eta_z^\theta([0,\tau])$ and the same boundary data as before on $\partial D$.

\subsection{Imaginary geometry lemmas}
The following lemma will allow us to compare flow lines generated from two instances of a GFF with different boundary values. We will only give a sketch of the proof, since the proof proceeds similarly as \cite[Lemma 5.4]{lqg-tbm2} (see also \cite[Remark 3.5]{IG1} for a related result).
\begin{lemma}\label{lem:RN}
	Let $h_1$ and $h_2$ be two Dirichlet GFF on $\C\setminus \D$ modulo a global additive multiple of $2\pi\chi$, such that\footnote{The maximal and minimal values of the field at $\partial \D$ are not well-defined, since the field is defined modulo $2\pi\chi$, but differences such that $h_i(x)-h_i(y)$ for $x,y\in\partial\D$ are well-defined. } $\sup_{x,y\in\partial \D}|h_i(x)-h_i(y)|<M$ for $i=1,2$ and some $M>0$. Let $U\subset\C\setminus \D$ be a domain bounded away from $\D$ and $\infty$. Then there exists a constant $C>0$ only depending on $M$ and $U$, such that $$\mu_1(E)\le C\mu_2(E)^{1/2},$$
	where $\mu_1$ and $\mu_2$ are the probability measures associated with $h_1$ and $h_2$, respectively, and $E$ is an arbitrary event in the $\sigma$-algebra of $h_1|_U$ (equivalently, $h_2|_U$).
\end{lemma}
\begin{proof}
	It is sufficient to prove the lemma under the assumption that $h_i$ is a Dirichlet GFF satisfying $\sup_{z\in\partial\D} |h_i(z)|<M$ for $i=1,2$, since the field $h_i$ in the statement of the lemma is a Dirichlet GFF, viewed modulo a global additive multiple of $2\pi\chi$. Let $g$ be the harmonic function in $\C\setminus \D$ which is constant at $\infty$, and has Dirichlet boundary values $(h_1-h_2)|_{\partial\D}$ on $\partial\D$. Then $h_1\eqD h_2+g$. Let $(\cdot,\cdot)_\nabla$ denote the Dirichlet inner product, which is defined by $(f_1,f_2)_\nabla:=(2\pi)^{-1}\int \nabla f_1\cdot\nabla f_2\,$ for smooth functions $f_1$ and $f_2$. It is explained in \cite[Lemma 5.4]{lqg-tbm2} that the Radon-Nikodym derivative of $h_1|_U$ with respect to $h_2|_U$ is given by
	\eqbn
	\exp( (h_2|_U,g|_U)_\nabla - \|g|_U\|_\nabla^2/2 ).
	\eqen
	For any event $E$ as in the statement of the lemma, the Cauchy-Schwarz inequality implies that
	\eqbn
	\mu_1(E) \leq \E[ \exp( 2(h_2|_U,g|_U)_\nabla - \|g|_U\|_\nabla^2 ) ]^{1/2}\cdot 
	\mu_2(E)^{1/2}.
	\eqen
	To conclude the proof, it is sufficient to show that the expected value on the right side is bounded by some constant only depending on $U$ and $M$. It is sufficient to show that $\|g|_U\|_\nabla\leq C$ for some constant $C$ satisfying these properties, since $(h_2|_U,g|_U)_\nabla$ is a normal random variable with variance $ \|g|_U\|_\nabla$ and expectation bounded in terms of $U$ and $M$. The result $\|g|_U\|_\nabla\leq C$ follows by standard regularity estimates for harmonic functions (see e.g.\ \cite[Chapter2, Theorem 7]{evans-pde}), which say that $|\nabla g|\leq C'$ for some $C'$ only depending on $U$ and $M$.
\end{proof}
The following basic lemma will be used later to deduce triviality of certain $\sigma$-algebras associated with whole-plane space-filling SLE. We remark that alternative arguments to prove similar results for other variants of the GFF can be found in e.g.\ \cite[Section 3.1]{IG1} and \cite[Lemma 8.2]{wedges}.
\begin{lemma}\label{lem-gfftrivial}
	Let $h$ be a whole-plane GFF modulo $2\pi\chi$.	For each $\delta>0$ let $\mcl F_{\delta}$ (resp., $\mcl G_R$)  be the $\sigma$-algebra generated by the restriction of $h$ to $B_\delta(0)$ (resp., $\C\setminus B_R(0)$). Then $\cap_{\delta>0} \mcl F_{\delta}$ and $\cap_{R>0} \mcl G_{R}$ are trivial. 
\end{lemma}
\begin{proof}
	Let $\wh h$ be a whole-plane GFF such that the average of $\wh h$ about $\partial\D$ is equal to zero, and let $U$ be an independent uniform random variable in $[0,2\pi\chi]$. Then $h$ is equal in law to $\wh h+U$ modulo $2\pi\chi$. The field $\wh h$ is invariant in law under the map $z\mapsto z^{-1}$, which implies that the same property holds for $h$. In order to conclude the proof it the lemma, it is therefore sufficient to show that $\cap_{R>0} \mcl G_{R}$ is trivial. Write $h=h^0+h^\dagger$, where $h^0$ is a radially symmetric function modulo $2\pi\chi$, and $h^\dagger$ is a distribution which has mean zero on any circle around the origin. Tail triviality of $h^0$ follows by using that $(h^0(e^{-t}))_{t\in\R}\eqD (B_t+U+2\pi\chi\Z)_{t\in\R}$, where $B$ is a standard two-sided Brownian motion. For $n\in\N$ let $\alpha_n$ be independent standard normal random variables, and let $(f_n)_{n\in\N}$ be an orthonormal basis for the Dirichlet inner product for the set of smooth compactly supported functions in $\C$ with mean zero. A whole-plane GFF $\wt h$ modulo a global additive constant can be written in the form $\sum_n \alpha_n f_n$, which implies that if $\wt{\mcl G}_R$ is the $\sigma$-algebra generated by the restriction of $\wt h$ to $\C\setminus B_R(0)$, then $\cap_{R>0}\wt{\mcl G}_R$ is trivial. Writing $\wt h=\wt h^0+\wt h^\dagger$ as above, it follows that the tail of $\wt h^\dagger$ is trivial. Since $\wt h^\dagger\eqD h^\dagger$, the tail of $h^\dagger$ is also trivial.
\end{proof}

\subsection{Space-filling SLE$_\kappa$}
For $\kappa>4$, whole-plane space-filling SLE$_{\kappa}$ is a space-filling curve in $\C$ which starts and ends at $\infty$. It is closely related to regular SLE$_{\kappa}$ by the following informal descriptions. For $\kappa\geq 8$ the law of a whole-plane space-filling SLE$_{\kappa}$ can be obtained by considering a regular chordal or radial SLE$_{\kappa}$ in any domain $D$, fixing some point $z\in D$ independent of the SLE, and ``zooming in'' near $z$. For $\kappa\in(4,8)$ we may define a chordal or radial space-filling SLE$_\kappa$ by considering a regular chordal or radial SLE$_\kappa$, and  filling in the created bubbles by independent space-filling loops. As above, we obtain whole-plane space-filling SLE$_{\kappa}$ by considering the local behavior of the chordal or radial space-filling curve near some fixed point.

Whole-plane space-filling SLE$_{\kappa}$ for all $\kappa>4$ was first constructed by using imaginary geometry with parameters as in \eqref{eq92}, see \cite{IG4} and \cite[Footnote 9]{wedges}. For any fixed $z_1,z_2\in\C$, the two flow lines $\eta_{z_1}^L$ and $\eta_{z_2}^L$ will eventually merge, and before this happens the curves will a.s.\ never cross each other. Therefore the set of flow lines $\eta_z^L$ for all $z\in\Q^2$ form a tree in $\C$ which is rooted at infinity, such that two branches in the tree never cross each other. The whole-plane space-filling SLE$_{\kappa}$ is defined to be the curve which traces this tree. More precisely, first define a total ordering on all points of $\Q^2$ by saying that $z_1$ comes before $z_2$ if $\eta_{z_1}^L$ merges into $\eta_{z_2}^L$ on the left side. A separate argument (see \cite[Section 4.3]{IG4})  shows that there is a well-defined continuous space-filling curve in $\C$ which visits the points of $\Q^2$ according to this order, and we define $\eta$ to be this curve.

\begin{lemma}
A whole-plane space-filling SLE$_\kappa$ $\eta$ parametrized by Lebesgue measure has stationary increments.
	\label{thm:ergodic}
\end{lemma}	
\begin{proof}
	We want to show that for any fixed $t\in \R$ we have $\eta\eqd \eta(\cdot+t)-\eta(t)$. The proof will proceed similarly as \cite[Lemma 9.3]{wedges}, where an analogous result for quantum parametrization of $\eta$ was shown. For any $z\in\C$, let $\tau_z:=\inf\{ t\in\R\,:\,\eta(t)=z \}$ be the first time at which $\eta$ hits $z$. For any fixed $R>0$, let $z_0$ be sampled uniformly at random from Lebesgue measure on $B_R(0)$, independently of $\eta$. By translation invariance in law of the GFF, and by independence of $z_0$ and $\eta$, we have $\eta(\cdot+\tau_{z_0})-z_0\eqd \eta$, which implies that
	\eqb
	\big( \eta,\eta(\cdot+t)-\eta(t )\big) \eqD
	\big(\eta(\cdot+\tau_{z_0})-z_0, \eta(\cdot+\tau_{z_0}+t)-\eta(\tau_{z_0}+t )\big).
	\label{eq91}
	\eqe  
	When $R\rta\infty$, the total variation distance between the laws of $z_0$ and $\eta(\tau_{z_0}+t)$, hence  $\tau_{z_0}$ and $\tau_{z_0}+t$, converges to 0. Therefore the total variation distance between the laws of $\eta(\cdot+\tau_{z_0})-z_0$ and $\eta(\cdot+\tau_{z_0}+t)-\eta(\tau_{z_0}+t )$ converges to 0. Since the laws of the two elements on the right side of \eqref{eq91} are arbitrarily close in total variation distance as $R\rta\infty$, we see that the two elements on the left side of \eqref{eq91} are equal in law. This implies the desired stationarity result.
\end{proof}

\section{Convergence of discrete contour function for $\kappa=8$}\label{sec:tight}
In this section we prove Theorem \ref{thm:convergence}. The main inputs to the proof are Theorem \ref{thm:lawler-viklund-wholeplane}, a chordal version of Proposition \ref{prop:lsw}, and Lemma \ref{lem:metric-compare}, which are results proved by other authors in \cite{LV-radial}, \cite{lsw-lerw-ust}, and \cite{barlow-sub}, respectively.

First we define a metric $\rho$ on the space of paths in $\C$.
For $i=1,2$ let $I_i\subseteq \R$ be an interval, and let $\gamma^i:I_i\to\C$ be a continuous function, i.e., $\gamma^i$ is a curve in $\C$. Then the distance $\rho(\gamma^1,\gamma^2)$ between $\gamma^1$ and $\gamma^2$ is given by 
\eqb
\begin{split}
	\rho(\gamma^1,\gamma^2) =&\, \inf \sum_{k=0}^{\infty} \min\left\{2^{-k};\, \sup_{t\in I_1\cap[-2^k,2^k]}|\alpha(t)-t| 
	+ \sup_{t\in I_1\cap[-2^k,2^k]} |\gamma^1(t)-\gamma^2(\alpha(t))| \right\},
\end{split}
\label{eq:rho}
\eqe
where the infimum is over all increasing homeomorphisms $\alpha:I_1\to I_2$. The following result is proved in \cite{LV-radial}.
\begin{theorem}[Lawler-Viklund'17]
	There is a universal constant $\check c>0$ such that for all $\ep>0$ and simply connected domains $D$ containing the origin with analytic boundary, there exists a $\delta_0\in(0,1]$ satisfying the following.  For each $\delta\in(0,\delta_0]$ consider a simple random walk on $\delta\Z^2$ started at 0 and run until hitting $\partial D$, and let $\ul\eta^\delta$ be the loop-erasure of the random walk. We view $\ul\eta^\delta$ as a continuous curve parametrized such that each edge is traversed in time $\check c\delta^{5/4}$. Let $\ul\eta$ be a radial SLE$_2$ in $D$ towards $0$, started from a point on $\partial D$ sampled from harmonic measure, and let $\ul\eta$ be parametrized by $5/4$-dimensional Minkowski content. This parametrization of $\ul\eta$ is well-defined, and there is a coupling of $\ul\eta^\delta$ and $\ul\eta$, such that 
	\eqbn
	\P[ \rho(\ul\eta,\ul\eta^\delta)>\ep ]<\ep.
	\eqen
	\label{thm:lawler-viklund-wholeplane} 
\end{theorem}
In the remainder of the section we let $\eta$ and $\eta^\delta$ for $\delta\in(0,1]$ be as in the statement of Theorem \ref{thm:convergence}, i.e., $\eta$ is a whole-plane space-filling SLE$_8$ parametrized by Lebesgue measure, and $\eta^\delta$ is the Peano curve of a uniform spanning tree on $\delta\Z^2$. A chordal version of the following proposition was proved in \cite{lsw-lerw-ust}.
\begin{proposition}
	For any $\ep>0$ we can find a $\delta_0>0$ such that for any $\delta<\delta_0$ there is a coupling of $\eta^\delta$ and $\eta$ satisfying $\P[ \rho(\eta,\eta^\delta)>\ep ]<\ep$.
	\label{prop:lsw}
\end{proposition}

We will first argue joint convergence of the uniform spanning tree and its dual in the topology introduced by Schramm in \cite{schramm0}. For any compact topological space $X$ let $\mcl H(X)$ be the set of compact subsets of $X$ equipped with the Hausdorff topology. Letting $\wh{\C}$ denote the Riemann sphere, define the topological space $\mcl O\mcl S$ by $\mcl O\mcl S=\mcl H(\wh{\C}\times\wh{\C}\times \mcl H(\wh{\C}))$. A spanning tree on $\delta\Z^2$ for some $\delta\in(0,1]$ can be represented by an element $\mcl T^\delta$ in $\mcl O\mcl S$ by saying that $(a,b,K)\in\mcl T^\delta$ iff $K$ is a simple path from $a\in\delta\Z^2$ to $b\in\delta\Z^2$ in the spanning tree. We let $\wt{\cT}^\delta$ denote represent the dual tree, and we denote the continuum analogues by $\cT$ and $\wt{\cT}$, respectively.
\begin{lemma}
	The pair $(\cT^\delta,\wt{\cT}^\delta)$ converge jointly to $(\cT,\wt{\cT})$ in $\cO\cS\times\cO\cS$
\end{lemma}
\begin{proof}
		 Tightness of $\cT^\delta$ is immediate since $\mcl O\mcl S$ is compact. For a UST on $\delta\Z^2$ and any finite collection of points $z_1,\dots,z_k\in\C$, let $\mcl T^\delta_{\bf z}$ for ${\bf z}=(z_1,\dots,z_k)$ be the element in $\mcl O\mcl S$ corresponding to the branches in the tree connecting $z_1,\dots,z_k$ (or the nearest lattice approximations of these points) to each other and to $\infty$. We define $\mcl T^\delta_{\bf z}$ similarly if ${\bf z}=(z_1,z_2,\dots)$ is countably infinite. An instance of a whole-plane space-filling SLE$_8$ $\eta$ in $\C$ gives elements $\mcl T$ and $\mcl T_{\bf z}$ in $\mcl O\mcl S$ by letting the branch or branches from each $z\in\C$ to $\infty$, be the left frontier of $\eta((-\infty,t])$ for each time $t$ satisfying $\eta(t)=z$. For each fixed $z$ there is a.s.\ only one such branch, and this branch is given by the flow line $\eta_z^L$ defined in Section \ref{sec:ig}. For any $R>1$, let $\mcl T^R,\mcl T^R_{\bf z},\mcl T^{R,\delta},\mcl T^{R,\delta}_{\bf z}$ be defined similarly, but for a chordal SLE$_8$ in $B_R(0)$ from $-Ri$ to $Ri$ (in the continuum case) or a UST on $B_R(0)\cap(\delta\Z^2)$ with half wired and half free boundary conditions (in the discrete case). By \cite{lsw-lerw-ust}  we know that $\mcl T^{R,\delta}_{\bf z}$ converges to $\mcl T^R_{\bf z}$ in $\mcl O\mcl S$. 
	
	By \cite[Corollary 4.5]{masson-lerwgrowth} the total variation distance between the laws of $\mcl T^{R,\delta}_{z_1}$ and $\mcl T^\delta_{z_1}$ for fixed $z_1\in\C$ goes to zero as $R\rta\infty$, uniformly in $\delta$. By Wilson's algorithm \cite{wilson-algorithm}, we get further that the total variation distance between the laws of $\mcl T^{R,\delta}_{\bf z}$ and $\mcl T^\delta_{\bf z}$ for any finite tuple $\bf z$ goes to zero as $R\rta\infty$, again uniformly in $\delta$.
	The total variation distance between the laws of $\mcl T^{R}_{\bf z}$ and $\mcl T_{\bf z}$ goes to zero as $R\rta\infty$, since this property holds for the Gaussian free fields from which the chordal and whole-plane, respectively, SLE$_8$'s were generated \cite[Proposition 2.11]{IG4}. We conclude that $\cT^\delta_{\bf z}$ converges in law to $\cT_{\bf z}$ in $\mcl O\mcl S$ for any finite tuple $\bf z$. By symmetry the same result holds for the dual, i.e., $\wt{\cT}^\delta_{\bf z}$ converges in law to $\wt{\cT}_{\bf z}$. 
	
	Let $(\cT',\wt{\cT}')$ be some subsequential scaling limit of the pair $(\cT^\delta,\wt{\cT}^\delta)$ in $\cO\cS\times\cO\cS$. We want to show that $(\cT',\wt{\cT}')\eqd (\cT,\wt{\cT})$. Let $\bf z$ be some enumeration of the rationals. By the convergence result for finite skeletons we see that $\cT'_{\bf z}\eqd \cT_{\bf z}$ and $\wt{\cT}'_{\bf z}\eqd \wt{\cT}_{\bf z}$. By \cite[Theorem 10.7]{schramm0} the trunk of $\cT'$ and the trunk of $\wt{\cT}'$ are disjoint. This gives that $\cT'_{\bf z}$ (resp.\ $\wt{\cT}'_{\bf z}$) uniquely determines $\wt{\cT}'$ (resp.\ $\cT'$), since the trunk of the trees are dense. Therefore $(\cT'_{\bf z},\wt{\cT}')\eqd (\cT_{\bf z},\wt{\cT})$ and $(\cT',\wt{\cT}'_{\bf z})\eqd (\cT,\wt{\cT}_{\bf z})$, which implies further that $(\cT',\wt{\cT}')\eqd (\cT,\wt{\cT})$. 
\end{proof}

\begin{proof}[Proof of Proposition \ref{prop:lsw}]
	Consider a coupling such that $(\cT^\delta,\wt{\cT}^\delta)$ converges to $(\cT,\wt{\cT})$ a.s.\ in $\cO\cS\times\cO\cS$. By the construction of space-filling SLE from imaginary geometry, the pair $(\cT,\wt{\cT})$ uniquely determines a space-filling curve $\eta$ with the law of a whole-plane space-filling SLE$_\kappa$. For any $\ep>0$ consider the flow lines $\eta^{L}_z$ and $\eta^{R}_z$ for $z\in \ep \Z^2$. The complement of these flow lines is a collection of open domains which we call continuum pockets, such that each pocket is enclosed by flow lines $\eta_{z_1}^L,\eta_{z_1}^R,\eta_{z_2}^L,\eta_{z_2}^R$ for $z_1,z_2\in\ep\Z^2$. Since the double points of SLE$_8$ have zero Lebesgue measure, for each fixed $z\in\C$ and any $\ep_1>0$, there is a.s.\ a  random  $\delta_1>0$, such that for all $w\in B_{\delta_1}(z)$, the flow line $\eta_w^{L}$ merges into $\eta_z^L$ before leaving $B_{\ep_1}(w)$. This implies that for any $z_1,z_2\in\ep\Z^2$, the Hausdorff distance between $\eta_{z_i}^{L,\delta}$ and $\eta_{z_i}^{L}$ restricted to any compact set, converges a.s.\ to zero, and, since $(\cT^\delta,\wt{\cT}^\delta)\rta(\cT,\wt{\cT})$ a.s., the point at which $\eta_{z_1}^{L,\delta}$ and $\eta_{z_2}^{L,\delta}$ merge, converges a.s.\ to the point at which $\eta_{z_1}^{L}$ and $\eta_{z_2}^{L}$ merge. It follows that a continuum pocket enclosed by flow lines $\eta_{z_1}^L,\eta_{z_1}^R,\eta_{z_2}^L,\eta_{z_2}^R$ for $z_1,z_2\in\ep\Z^2$, is a.s.\ the limit for the Hausdorff distance of a discrete pocket enclosed by the corresponding discrete flow lines. The Peano curves $\eta^\delta$ and $\eta$ visit the pockets in an order corresponding to tracing the interface of the primal tree and the dual tree, and the order in which the pockets are visited, converges a.s.\ as $\delta\rta 0$, if we only consider the pockets restricted to some compact set. Therefore, for any fixed $T>0$ and with $p_\ep$ being the maximal diameter of the continuum pockets visited by $\eta$ during $[-T,T]$, we have $|\eta(t)-\eta^\delta(t)|<10p_\ep$ a.s.\ for all sufficiently small $\delta$ and all $t\in[-T,T]$. Since $\lim_{\ep\rta 0}p_\ep= 0$, we see upon decreasing  $\ep$ that $\lim_{\delta\rta 0}\rho(\eta,\eta^\delta)=0$.	
\end{proof}

We recall the following result from \cite{barlow-sub}. For any $\delta\in(0,1]$ and a set of edges $A$ of the square grid $\delta\Z^2$, we define $\frk m^\delta(A):=|A|\check c\delta^{5/4}$, with $\check c$ as in Theorem \ref{thm:lawler-viklund-wholeplane} and $|A|$ denoting the number of elements in $A$. For $\delta\in(0,1]$, a UST on $\delta\Z^2$, and $z\in\C$, let $\eta^{L,\delta}_z$ be the path in the UST  from the nearest lattice approximation of $z$ to $\infty$.
\begin{lemma}[Proposition 2.8, \cite{barlow-sub}] \label{lem:metric-compare}
	There exist universal constants $c_1,c_2,c_3,\lambda_0>0$ such that the following is true for any $\delta\in(0,1]$. Given $r\ge \delta$ and $\lambda\geq\lambda_0$, let $R=re^{c_1\lambda^{1/2}} $. Let $A(r,\lambda)$ be the event that  
	for all $x,y\in B_R(0)\cap (\delta \Z^2)$ such that $\op{diam}( \eta_x^{L,\delta}\Delta \eta_y^{L,\delta} ) \le r$, we have $\frk{m}^\delta( \eta_x^{L,\delta}\Delta \eta_y^{L,\delta} ) \le \lambda r^{5/4}$, where $\Delta$ denotes symmetric difference.
	For every $r\ge \delta$ and $\lambda\geq\lambda_0$ we have $\P[A(r,\lambda)^c]\le c_3\exp\{-c_2\lambda^{1/2}\}$.
\end{lemma} 

Next we prove tightness of the rescaled version $Z^\delta$ of $\frk Z$. Recall the definition of $Z^\delta=(L^\delta,R^\delta)$  above the statement of Theorem \ref{thm:convergence}.
\begin{proposition}\label{prop:tight}
	The contour functions $Z^\delta$ for $\delta\in(0,1]$ are tight for the topology of uniform convergence on compact sets.
\end{proposition}

\begin{proof}
	By scale invariance, it is sufficient to show that $Z^{\delta}|_{[0,1]}$ is tight, and by symmetry in $L^\delta$ and $R^\delta$ it is sufficient to prove that $L^\delta$ is tight. For any $\delta>0$, let $w^\delta:[0,1]\to[0,\infty)$ be the minimal increasing modulus of continuity of $L^\delta|_{[0,1]}$, i.e., it is the minimal increasing function such that for any $t,s\in[0,1]$ we have
	\eqbn
	|L^\delta_t-L^\delta_s|\leq w^\delta(|t-s|).
	\eqen 
	By \cite[Lemma 2.1]{prohorov-tight} and since $L^\delta_0=0$ for all $\delta\in(0,1)$, tightness of $L^{\delta}|_{[0,1]}$ follows if we can prove that for any $\ep>0$ there exists a $\rho>0$ such that $\P[w^\delta(\rho)>\eps]<\eps$ for all $\delta\in(0,1]$. 
	Since $L^\delta$ is Lipschitz continuous with constant $\check c^{-1}\delta^{-5/8}$, it is sufficient to show that this holds for small $\delta$, i.e., it is sufficient to show that for all $\ep>0$ there exists $\rho,\delta_0>0$ such that $\P[w^\delta(\rho)>\eps]<\eps$ for all $\delta\in(0,\delta_0]$. Choose $\lambda>0$ sufficiently large and $r>0$ sufficiently small such that, in the notation of Lemma \ref{lem:metric-compare} and for all sufficiently small $\delta\in(0,1]$, we have $\lambda r^{5/4}<\ep$, $\P[\op{diam}(\eta([0,1]))>R/2]<\ep/4$ and $\P[A(r,\lambda)^c]<\ep/3$. By Proposition \ref{prop:lsw}, there exists a $\delta_0>0$ such that for $\delta<\delta_0$ we have $\P[\op{diam}(\eta^\delta([0,1]))>R]<\ep/3$. Since $\eta$ is continuous a.s., we may choose $\rho>0$ sufficiently small such that 
	\eqbn
	\P\left[ \sup_{s,t\in[0,1],0<t-s<2\rho} \op{diam} \eta([s,t]) > r/2 \right]<\ep/4.
	\eqen
	Applying Proposition \ref{prop:lsw} again and decreasing $\delta_0>0$ if necessary, the following holds for any $\delta<\delta_0$
	\eqbn
	\P\left[ \sup_{s,t\in[0,1],0<t-s<\rho} \op{diam} \eta^\delta([s,t]) > r \right]<\ep/3.
	\eqen
	Combining the above results, we have shown that with probability at least $1-\ep$, the following event occurs
	\eqb
	\{ \op{diam}(\eta^\delta([0,1]))\leq R \} \cap
	A(r,\lambda) \cap
	\left\{ \sup_{s,t\in[0,1],0<t-s<\rho} \op{diam} \eta^\delta([s,t]) \leq r \right\}.
	\label{eq93}
	\eqe
	On the event \eqref{eq93}, for all $s,t\in[0,1]$ such that $0<t-s<\rho$, we have $\op{diam}(\eta^\delta([s,t]))\leq r$, which implies by occurrence of $A(r,\lambda)$ that $\frk{m}^\delta( \eta_{\eta(t)}^{L,\delta}\Delta \eta_{\eta(s)}^{L,\delta} ) \le \lambda r^{5/4}<\ep$. Therefore  $w^\delta(\rho)<\ep$, and the lemma follows.
\end{proof}

\begin{proposition}
	The pair
	$(\eta^\delta,Z^\delta)$ converges weakly to $(\eta,Z)$.
\end{proposition}
\begin{proof}	
	By Propositions \ref{prop:lsw} and \ref{prop:tight}, the pair $(\eta^\delta,Z^\delta)$ converges subsequentially in law to some limiting random variable $(\eta,\wt Z)$, where $\eta$ has the law of an SLE$_8$, and $\wt Z=(\wt L,\wt R)$ is continuous. Considering a coupling for different $\delta$ where this subsequential convergence holds a.s., we need to prove that $\wt Z=Z$ a.s., where $Z$ is as in the statement of the proposition. For any $z\in\C$ let $\tau(z):=\inf\{t\in\R\,:\,\eta(t)=z \}$ be the first time at which $\eta$ hits $z$. We observed in the proof of Proposition \ref{prop:lsw} that $\eta^\delta$ converges jointly with the finite skeletons $\mcl T_{{\bf z}}^\delta$ and $\wt{\mcl T}_{{\bf z}}^\delta$. By Theorem \ref{thm:lawler-viklund-wholeplane} and since the natural parametrization of $\ul\eta^\delta$ (resp.\ $\ul\eta$) is determined by the unparametrized curve, we have joint convergence in law of $\eta^\delta$ and the branches of $\mcl T_{{\bf z}}^\delta$ and $\wt{\mcl T}_{{\bf z}}^\delta$ viewed as \emph{parametrized} curves.
	This implies that for any fixed $z\in\C$, $(\eta^\delta,L^\delta_{\tau(z)})$ converges in distribution to $(\eta,L_{\tau(z)})$. Since $\wt L$ is a.s.\ continuous and $\{ \tau(z)\,:\,z\in\Q^2 \}$ is dense in $\R$, we see that $\wt L=L$ a.s. We have $\wt R=R$ a.s.\ by a similar argument, which completes the proof. \qedhere
	
\end{proof}

\section{Existence and properties of the contour functions}\label{sec:contour}
In this section we will prove Theorem \ref{thm:Z}, which says that the contour functions $Z=(L,R)$ for all $\kappa>4$ are well-defined and satisfy certain basic properties.

We will prove that $Z$ is well-defined as a continuous function by using the Kolmogorov-Chentsov theorem, and we therefore need a moment bound for the increments of $Z$. We will obtain a moment bound by drawing the space-filling SLE$_\kappa$ $\eta$ on top of a $4/\sqrt{\kappa}$-LQG surface, and using that the Minkowski content of the SLE frontier is given by the expected quantum length of the frontier, up to multiplication by a function depending on local properties of the field.

Let $\ul\kappa\in(0,4)$, $\gamma=\sqrt{\ul\kappa}$ and $D\subseteq\C$, and let $h$ be some GFF-like field on $D$. Recall from the introduction that Liouville quantum gravity (LQG) with parameter $\gamma$ is a random surface associated with $h$. In particular, the field $h$ induces a random area measure $\mu_h$ on $D$ which may be written heuristically in the form $e^{\gamma h}\,dz$, where $dz$ is Lebesgue measure \cite{shef-kpz,rhodes-vargas-review}.

The field $h$ associated with a $\gamma$-LQG surface also induces a length measure along certain curves, e.g.\ along $\partial D$ or SLE$_{\ul\kappa}$ curves in $D$. For an SLE$_{\ul\kappa}$ or SLE$_{\ul\kappa}(\ul\rho)$ curve $\ul\eta$ in $D$ there are two natural ways to define such a $\gamma$-LQG length measure. The first approach is to define a measure $\nu_h$ on $\ul\eta$ by considering the quantum boundary length measure as defined in e.g.\ \cite{shef-kpz,Zipper}. We consider a conformal map $\psi:U\to\BB H$, where $U$ is some domain on one ``side'' of $\ul\eta$, such that $\psi$ straightens $\ul\eta$. Consider the $\gamma$-LQG boundary measure on $\R$ which we get when applying the coordinate change formula for quantum surfaces to $h$ and $\psi$. Let $\nu_h$ be the pullback under $\psi$ of this quantum measure on $\R$. Note that we may view $\nu_h$ as a measure on $\C$ supported on $\ul\eta$.

The second approach is to define a $\gamma$-LQG measure $\sigma_h$ with (roughly speaking) the Minkowski content of $\ul\eta$ as base measure. Recalling that $\frk m$ denotes $(1+\ul\kappa/8)$-dimensional Minkowski content and considering some arbitrary strictly monotone parametrization $\ul\eta:\R_+\to\C$ of $\ul\eta$ satisfying $\ul\eta(0)=0$, we first define the measure $\frk m_{\ul\eta}$ on $\C$ by
\eqb
\frk m_{\ul\eta}(U) = \lim_{\ep\rta 0} \frk m( U\cap \ul\eta([\ep,\infty) )).
\label{eq:contour12}
\eqe
\begin{remark}\label{rmk:origin}
We define the measure $\frk m_{\ul\eta}$ as a limit, rather than considering $\frk m( U\cap \ul\eta([0,\infty) ))$, since it is not known that the Minkowski content of $\ul\eta$, which is a whole plane $\SLE_{\ul\kappa}(2-\ul\kappa)$, is well-defined near the origin. This assertion holds for whole plane $\SLE_{\ul \kappa}$ (see \cite[Lemma 3.1]{LV-radial} and \cite{whole-plane-zhan}), which means for $\ul\kappa=2$ we can remove the cutoff in \eqref{eq:contour12}.
\end{remark}
For $z\in D$ and $\ep>0$ such that $B_\ep(z)\subset D$, we let $h_\ep(z)$ denote the average of $h$ around the circle $\partial B_\ep(z)$, see \cite{shef-kpz}. The measure $\sigma_h$ is defined by the following limit for any open set $U$ 
\eqb
\sigma_h(U) = \lim_{\ep\rta 0} \int_U e^{\gamma h_\ep(z)/2} \ep^{\gamma^2/8}\,d\frk m_{\ul\eta}(z).
\label{eq:contour4}
\eqe
When $h$ is a centered Gaussian field, the convergence holds  in $L^1$ for any bounded Borel set $U$ \cite[Theorem 1.1]{berestycki-elem}. See  \cite{berestycki-elem} and \cite{benoistnotes} for further details about $\sigma_h$. 
\begin{lemma}
	Let $h$ be a whole-plane GFF such that the average of the field over the unit circle is zero. Let $\ul\eta$ be an independent whole-plane SLE$_{\ul\kappa}(2-\ul\kappa)$, and let $\sigma_h$ and $\nu_h$ be as above. Then there exists a deterministic constant $c>0$ so that  $\sigma_h=c\nu_h$ a.s.
	
	 Furthermore, for any $a>0$ there is a $C>0$ such that if $\mu_h$ is the $\gamma$-LQG area measure associated with $h$ and  $U\subset \D$, then we have $\frk m_{\ul\eta}(U)\leq C\E[ \sigma_h(U)\1_{\mu_h(\D)<a} \,|\,\ul\eta]$.
\end{lemma}
\begin{proof}
We first prove the following claim \eqref{eq:benoist}. Let $\gamma=\sqrt{\ul\kappa}$ and $\wt h$ be a free boundary GFF plus $(\gamma -\frac{2}{\gamma})\log|z|^{-1}$ in $\BB H$, with some arbitrary choice of additive constant, and let $\wt{\ul\eta}$ be an independent SLE$_{\ul\kappa}$ in $\BB H$ from 0 to $\infty$. Define $\frk m_{\wt\eta}$, $\nu_{\wt h}$ and $\sigma_{\wt h}$ similarly as in \eqref{eq:contour12} and \eqref{eq:contour4}. Then there is a constant $c>0$ such that 
\begin{equation}\label{eq:benoist}
\sigma_{\wt h}=c\nu_{\wt h} \qquad \textrm{a.s.}
\end{equation}
 In fact \eqref{eq:benoist} is proved  in \cite[Proposition 3.3]{benoistnotes} if  $\wt h$ is replaced by a so-called $(\gamma -\frac{2}{\gamma})$-quantum wedge.
By the definition of quantum wedge and its relation to $\wt h$ (see e.g.\ \cite[Section 4.2]{wedges}), \eqref{eq:benoist} holds for $\wt h$. 
Indeed, the quantum wedge is defined up to a scaling of the complex plane and the free GFF is defined up to an additive constant. For any $R>0$, it is possible to choose the scaling for the quantum wedge and	the additive constant for the free GFF so that $\wt h$ agree with the  $(\gamma -\frac{2}{\gamma})$-quantum wedge on $B_R(0)$.
	
Both the measures $\nu_h$ and $\sigma_h$ are defined locally, in the sense that for an arbitrary monotone parametrization of $\ul\eta$ and any interval $I\subset\R_+$, the measure of $\ul\eta(I)$ depends only on $\ul\eta(I)$ and on $h$ restricted to some neighborhood of $\ul\eta(I)$. Bounded away from $\partial \D$, 0 and $\infty$, the field $h$ is absolutely continuous with respect to translations of the field $\wt h$, and $\ul\eta$ is locally absolutely continuous with respect to the curve $\wt{\ul\eta}$, in the sense that any interval $I\subset\R_+$ bounded away from 0 and $\infty$ can be written as a finite union of (random) intervals $I_i$, such that $\ul\eta|_{I_i}$ is absolutely continuous with respect to a segment of $\wt{\ul\eta}$.
	This implies that $\sigma_h=c\nu_h$ a.s., where $c$ is as in \eqref{eq:benoist}. This proves the first assertion.
	
	Now we prove the second assertion. By \eqref{eq:contour4}, for any $U\subset\D$,
	\eqbn
	\begin{split}
		\E[\sigma_h(U)\1_{\mu_h(\D)<a}\,|\,\ul\eta]
		= \lim_{\ep\rta 0} \E\left[ \int_U e^{\gamma h_\ep(z)/2}\ep^{\gamma^2/8} \,d\frk m_{\ul\eta}(z) \1_{\mu_h(\D)<a}\,|\,\ul\eta\right]. 
	\end{split}
	\eqen
	By independence of $h$ and $\ul\eta$, in order to conclude the proof, it is therefore sufficient to prove the existence of constants $c_1,\ep_0>0$, such that for any $z\in U$ and $\ep\in(0,\ep_0)$,
	\eqb
	\begin{split}
		\E\left[e^{\gamma h_\ep(z)/2}\ep^{\gamma^2/8} \1_{\mu_h(\D)<a}\right] \geq c_1.
	\end{split}
	\label{eq:contour9}
	\eqe
	For any $z\in\C$ let $h^z:= h(\cdot+z)-h_1(z)$, and observe that $h^z\eqd h$. For any field $\ol h$ and for fixed $c_2,c_3>0$ let $E(\ol h)$ be the event that $\mu_{\ol h}(B_{10}(0))\leq c_2$, and that $e^{\gamma\ol h_1(w)}\in[c_3^{-1},c_3]$ for all $w\in B_{10}(0)$. Choose $c_2,c_3$ such that $c_2c_3< a$ and $\P[E(h)]>0$, and observe that appropriate constants exist since $\mu_{h}(B_{10}(0))<\infty$ a.s., and $h_1(w)$ is a.s.\ continuous in $w$ \cite[Proposition 3.1]{shef-kpz}. For any $z,w\in \D$ and $\ep>1$, since $h=(h^{z-w})^{w-z}$ implies that $ h_\ep(z)=h_\ep^{z-w}(w)-h_1^{z-w}(-z+w)$, 
	\begin{align}
	\E\left[\ep^{\frac{\gamma^2}{8}} e^{\frac{\gamma}{2} h_\ep(z)} \1_{\mu_{h}(\D)\leq a}  \right]
	&= \E\left[e^{-\frac{\gamma}{2} h^{z-w}_1(-z+w)}\cdot \ep^{\frac{\gamma^2}{8}} e^{\frac{\gamma}{2} h_\ep^{z-w}(w)} \1_{e^{-\gamma  h^{z-w}_1(-z+w)}\mu_{ h^{z-w}}(B_1(-z+w))\leq a}\right]\nonumber\\
	&= \E\left[e^{-\frac{\gamma}{2} h_1(-z+w)}\cdot \ep^{\frac{\gamma^2}{8}} e^{\frac{\gamma}{2} h_\ep(w)} \1_{e^{-\gamma h_1(-z+w)}\mu_{ h}(B_1(-z+w))\leq a}\right]\nonumber\\
	&\geq \E\left[c_3^{-1} \ep^{\frac{\gamma^2}{8}} e^{\frac{\gamma}{2}  h_\ep(w)} \1_{E( h)}
	\right].
	\label{eq:contour8}
	\end{align}
	Let $\wt\mu_{h}$ be the $\gamma/2$-LQG area measure associated with $h$. Since the regularized measures $\ep^{\frac{\gamma^2}{8}} e^{\frac{\gamma}{2} h_\ep(w)}\,dw$ converge to $\wt\mu_{ h}$ in $L^1$,
	\eqbn
	\E[\wt\mu_{ h}(\D)\1_{E(h)}] = \lim_{\ep\rta\infty} \int_{\D} \E\left[\ep^{\frac{\gamma^2}{8}} e^{\frac{\gamma}{2} h_\ep(w)}\1_{E(h)}
	\right]\,dw.
	\eqen
	Since $\wt\mu_{ h}(\D)>0$ a.s., there are $\ep_0,c_4>0$ such that for any $\ep\in(0,\ep_0)$ we can find a $w\in\D$ such that $\E\big[\ep^{\frac{\gamma^2}{8}} e^{\frac{\gamma}{2} h_\ep(w)}\1_{E(h)}
	\big]>c_4$. For such $\ep$ and $w$ we get by insertion into \eqref{eq:contour8} that
	$\E\big[\ep^{\frac{\gamma^2}{8}} e^{\frac{\gamma}{2} h_\ep(z)} \1_{\mu_{h}(\D)\leq a}  \big] \geq c_3^{-1}c_4$, so \eqref{eq:contour9} holds with $c_1=c_3^{-1}c_4$.
\end{proof}

In the next few paragraphs we let $\kappa>4$, and we let $\eta$ be a whole-plane space-filling SLE$_\kappa$ from $\infty$ to $\infty$. We will prove existence of the boundary length process $Z$ at a fixed time, and prove a moment bound for $Z$. By symmetry in $L$ and $R$, it is sufficient to consider $L$.
Recalling that $L$ describes the evolution of the length of the left frontier of $\eta$, we see that for any fixed $t\in\R$ and some arbitrary parametrization of $\eta^L_{\eta(t)}$ and $\eta^L_0$, a.s.
\eqb
L_t = \lim_{\ep\rta 0} \Big(\frk m\big( \eta^L_{\eta(t)}((\ep,\infty)) \setminus \eta^L_{0} \big)
- \frk m\big( \eta^L_{0}((\ep,\infty)) \setminus \eta^L_{\eta(t)} \big)\Big).
\label{eq:contour5}
\eqe
Remark~\ref{rmk:origin} explains why we define $L$ as a limit, rather than considering the Minkowski content of the full frontier.
\begin{lemma}
	The random variable $L_1$ defined by \eqref{eq:contour5} is well-defined a.s., and
	$\E[\sup_{t\in[0,1]}L_t^N]<\infty$ for all $N\in\N$.
	\label{prop30}
\end{lemma}

\begin{proof}
	Existence of right side of \eqref{eq:contour5} before we take the limit $\ep\rta 0$ follows by existence of the Minkowski content of chordal SLE$_{16/\kappa}$ and local absolute continuity. 
	
	A $\gamma$-quantum cone is a particular kind of $\gamma$-LQG surface, which may be constructed by sampling a point from the $\gamma$-LQG area measure of some $\gamma$-LQG surface, and zooming in near the sampled point. See \cite[Section 4]{wedges} for the formal definition and basic properties of a $\gamma$-quantum cone. Let $h$ be the field associated with a $\gamma$-quantum cone on $\C$, embedded such that the average of $h$ about $\partial \D$ is equal to zero. Let $\wt Z=(\wt L_t,\wt R_t)_{t\in\R}$ describe the evolution of the quantum boundary length of $\eta$, corresponding to a time parameterization of $\eta$ by quantum area. By Theorem \ref{thm:wedges} from \cite{wedges}, $\wt Z$ has the law of a two-dimensional correlated Brownian motion. We may couple $h$ with a whole-plane GFF $\wt h$ with unit circle average 0, such that $h|_{\D}=(\wt h-\gamma\log|\cdot|)|_{\D}$. 
	
	For $r_2>r_1>0$ let $A(r_1,r_2):=\{ z\in\C\,:\,r_1<|z|<r_2 \}$. 
	Define $T_1:=\inf\{t\geq 0\,:\, \eta(t)=3/4 \}$ and $T_2:=\inf\{t\geq T_1\,:\, \eta(t)\not\in B_{1/4}(3/4) \}$.
	Let $\tau_r=\inf\{t:|\eta(t)|=r \}$ for $r>0$. We now show that 
 	\begin{equation}\label{eq:key}
 			\P\left[ \sup_{t\in[0,\tau_{1/4}]} |L_t|>M \right]=\P\left[ \sup_{t\in[T_1,T_2]} |L_t-L_{T_1}|>M \right]\;\textrm{decays super-polynomially}.
 	\end{equation}
	Note that we have equality of the two probabilities in \eqref{eq:key} by  invariance of whole-plane SLE$_\kappa$ under recentering at a deterministic point, which follows from the analogous property of the whole-plane Gaussian free field modulo $2\pi\chi$.
	For fixed $t\in\R$ let $\sigma_{h;t}$ denote the measure $\sigma_h$ defined by \eqref{eq:contour4} with $\frk m_{\eta^L_{\eta(t)}}$ as base measure.
	By Lemma \ref{prop30}, there is a $C>0$ such that for any $U\subset A(1/2,1)$ and $t\in\R$ (both chosen in a way measurable with respect to $\eta$),
	\eqbn
	\frk m_{\eta_0^L}( U ) \leq C \E[ \sigma_{h;0}( U ) \1_{\mu_h(A(1/2,1))<1}\,|\,\eta ],\qquad 
	\frk m_{\eta_{\eta(t)}^L}( U ) \leq C \E[ \sigma_{h;t}( U )\1_{\mu_h(A(1/2,1))<1}\,|\,\eta ],
	\eqen
	so with $\frk t(t) := \op{sign}(t)\cdot \mu_h\big(\eta([0\wedge t,0\vee t])\big)$ and for any $t\in\R$,
\eqbn
|L_t-L_{T_1}| \leq C\E\left[ \big| \wt L_{\frk t(t)} - \wt L_{\frk t(T_1)}\big|\1_{\mu_h(A(1/2,1))<1} \,\,\Big|\,\,\eta \right],
\eqen	
and further
	\eqb
	\begin{split}
	\sup_{t\in[T_1,T_2]} |L_t-L_{T_1}| 
	&\leq
	C\sup_{s,t\in[T_1,T_2] }\E\left[  \big| \wt L_{\frk t(t)} - \wt L_{\frk t(s)} \big|	\1_{\mu_h(A(1/2,1))<1} \,\,\Big|\,\,\eta \right]\\
	&\leq 
	C\E\left[ \sup_{s,t\in[T_1,T_2] } | \wt L_{\frk t(t)} - \wt L_{\frk t(s)}|	\1_{\mu_h(A(1/2,1))<1} \,\,\Big|\,\,\eta \right].
	\end{split}
	\label{eq94}
	\eqe
	By an application of Chebyshev's inequality it follows that
	\eqb
	\begin{split}
		&\P\left[ \sup_{t\in[T_1,T_2]} |L_t-L_{T_1}|>M \right]\leq \left(\frac{C}{M}\right)^N \E\left[\E\left[ \sup_{s,t\in[T_1,T_2] } | \wt L_{\frk t(s)} - \wt L_{\frk t(t)}|	\1_{\mu_h(A(1/2,1))<1} \,\,\Big|\,\,\eta \right]^N \right] \\
		&\qquad\preceq M^{-N} \E\left[
		\sup_{\substack{\text{$t\in[0,\mu_h(\eta([0,T_2]))]$,} \\ \text{$s\in[t,t+1]$}}}
		|\wt L_{s} 
		-\wt L_{t}|^N
		\right] \\
		&\qquad\preceq M^{-N} \int_{\R_+} \P\left[ 
		\sup_{\substack{\text{$t\in[0,\mu_h(\eta([0,T_2]))]$,} \\ \text{$s\in[t,t+1]$}}}
		|\wt L_{s} -\wt L_{t}|>a  \right]a^{N-1}\,da \\
		&\qquad \leq M^{-N} \int_{\R_+} \P\left[\mu_h(\eta([0,T_2]))>a^{K_1} \right] a^{N-1} + 
		\P\left[
		\sup_{\substack{\text{$t\in[0,a^{K_1}]$,} \\ \text{$s\in[t,t+1]$}}}
		|\wt L_t-\wt L_s|>a \right] a^{N-1} \,da,
	\end{split}
	\label{eq:contour11}
	\eqe
	where the implicit constant depends on $N$ and $\gamma$. We consider each term in the integrand on the right side separately. For any $K_2>0$,
	\eqbn
	\P\left[ \mu_h(\eta([0,T_2]))>a^{K_1} \right]
	\leq \P[\eta([0,T_2])\not\subset B_{a^{K_2}}(0) ] + 
	\P[ \mu_h( B_{a^{K_2}}(0)) > a^{K_1} ].
	\eqen
	By Proposition \ref{prop:powerlaw} first term on the right side is $\preceq a^{-(N+10)}$ for sufficiently large $K_2$. By conformal invariance of the GFF and  since there exists $p>0$ such that $\E[\mu_h(\D)^p]<\infty$  (see the argument of \cite[Lemma A.1]{rhodes-vargas-sphere} for a proof), the second term on the right side is $\preceq a^{-(N+10)}$ if we choose $K_1$ sufficiently large after fixing $K_2$. By a union bound and the Markov property of Brownian motion,
	\eqbn
	\P\left[ \sup_{t\in[0,a^{K_1}],s\in[t,t+1]} |\wt L_t-\wt L_s|>a \right]
	\preceq a^{K_1} \P\left[ \sup_{t\in[0,1],s\in[t,t+1]} |\wt L_t-\wt L_s|>a \right],
	\eqen
	which is $\preceq a^{-(N+10)}$ by tail estimates for Brownian motion. Inserting these estimates into \eqref{eq:contour11}, we get
	\eqb
	\P\left[ \sup_{t\in[T_1,T_2]} |L_t-L_{T_1}|>M \right]\preceq M^{-N},
	\label{eq:contour13}
	\eqe
	where the implicit constant depends on $N$ and $\kappa$. Equation \eqref{eq:contour13} combined with translation  invariance of $\eta$ concludes \eqref{eq:key}.
	
	By a union bound,
	\eqb
	\P\left[ \sup_{t\in[0,1]} |L_t| > M \right]
	\leq 
	\P\left[ \op{diam}\eta([0,1])>M^{0.01} \right]
	+ \P\left[ \sup_{t\in[0,\tau_{M^{0.01}}]} |L_t|>M \right].
	\label{eq:contour2}
	\eqe
	By \cite[Lemma 3.6]{kpz-bm}, the first term on the right side decays faster than any power of $M$. By \eqref{eq:key}, along with scale invariance of space-filling SLE, the second term on the right side decays faster than any power of $M$. It follows that $\E[\sup_{t\in[0,1]}L_t^N]<\infty$ for all $N\in\N$.
\end{proof}

\begin{proof}[Theorem \ref{thm:Z}, (i) and (ii)]
	By scale invariance and translation invariance of SLE and of the Minkowski content, and by Lemma \ref{prop30}, we have $\E[|L_t-L_s|^N]\preceq |t-s|^{N(1/2+1/\kappa)}$ for any $N\in\N$ and $t,s\in\R$, where the implicit constant depends on $\kappa$ and $N$. We get the exponent $1/2+1/\kappa$ by scale invariance of SLE, and since $d$-dimensional Minkowski content is multiplied by $r^{d/2}$ under the map $z\mapsto r^{1/2}z$ for some $r>0$. The same result holds for $R$ instead of $L$. A quantitative version of the Kolmogorov-Chentsov theorem (see e.g.\ \cite[Proposition 2.3]{lqg-tbm2}) now implies that there is a function $Z$ satisfying (i) and the scaling result of (ii), such that for any given $t\in\R$, $L_t$ is given by \eqref{eq:contour5} a.s., and that the same result holds with $R$ instead of $L$. The stationary and tail triviality results of (ii) follow from Lemmas \ref{lem-gfftrivial} and \ref{thm:ergodic}. The result \eqref{eq:intro1} follows from scale invariance and tail triviality.
\end{proof}
The lower bound for the Minkowski content of a whole-plane SLE$_2$ in the following lemma will be used in the proof of Theorem~\ref{thm:Z} (iii). A similar super-polynomial lower bound for whole plane $\SLE_{\ul \kappa}(2-\ul\kappa)$  with other $\ul\kappa<4$, which we expect to hold, would imply that Theorem \ref{thm:Z} (iii) also holds for $\kappa\neq 8$.
\begin{lemma}
	Let $\ul\eta$ be a whole-plane $\SLE_2$ from 0 to $\infty$ in $\C$ with some arbitrary strictly monotone parametrization, and define $\tau:=\inf\{ t\geq 0\,:\,\ul\eta(t)\not\in\D \}$. There is a  constant $c>0$ such that for all $M>0$, and with $\frk m_{\ul\eta}$ defined by \eqref{eq:contour12},
	\eqb
	\P[ \frk m_{\ul\eta}(\ul\eta([0,\tau]))>M ]<2\exp(-cM).	\label{eq:contour7}
	\eqe
	For any $\alpha<4/5$ there are constants $C,c'>0$ such that for all $M>0$,
	\eqb
	\P[ \frk m_{\ul\eta}(\ul\eta([0,\tau]) )<M^{-1} ]<C\exp(-c' M^\alpha).
	\label{eq:contour3}
	\eqe
	\label{prop10}
\end{lemma}
\begin{proof}
	We will only give a proof of \eqref{eq:contour3}, since \eqref{eq:contour7} is proved in the exact same way. For $\delta\in(0,1]$ let $\ul\eta^\delta$ be a LERW on $\delta\Z^2$ from 0 to $\infty$. Fix $\alpha<4/5$, and define $\tau^\delta:=\inf\{ t\geq 0\,:\,\ul\eta^\delta(t)\not\in B_{1/2}(0) \}$. By \cite{barlow-exp-lerw} there are constants $C,c'>0$ such that for all $\delta\in(0,1]$,
	\eqbn
	\P\left[ \frk m^\delta\left(\ul\eta^\delta([0,\tau^\delta]) \right)< 2M^{-1} \right]<\frac 12 C\exp(-c' M^\alpha),
	\eqen
	where $\frk m^\delta$ is as defined above the statement of Lemma \ref{lem:metric-compare}. By Proposition \ref{prop:lsw}, given any $M>2$, there is a $\delta_M>0$, such that for any $\delta<\delta_M$ there is a coupling of $\ul\eta^\delta$ and $\ul\eta$ satisfying
	\eqbn
	\P\left[ \rho(\ul\eta^\delta,\ul\eta)> \frac{1}{10M} \right]<\frac 12 C\exp(-c' M^\alpha).
	\eqen
	To conclude the proof of \eqref{eq:contour3}, it is therefore sufficient to show that 
	\eqbn
	\{\frk m_{\ul\eta}(\ul\eta([0,\tau]))<M^{-1}\} \subset \left\{ \frk m^\delta(\ul\eta^\delta([0,\tau^\delta]))< 2 M^{-1} \right\}
	\cup \left\{ \rho(\ul\eta^\delta,\ul\eta)> \frac{1}{10M}   \right\}.
	\eqen
	We will prove this result by contradiction, and we assume the event on the left side occurs, but neither of the two events on the right side occurs. Choose a homeomorphism $\alpha:\R_+\to\R_+$ such that the right side of \eqref{eq:rho} differs from $\rho(\ul\eta^\delta,\ul\eta)$ by less than $\frac{1}{10M}$. We obtain a contradiction by observing that $|\alpha(t)-t|$ (resp.\ $|\ul\eta^\delta(t)-\ul\eta(\alpha(t))|$) is larger than $\frac{1}{5M}$ for $t=\alpha^{-1}(\tau)$ if $\alpha^{-1}(\tau)>2M^{-1}$ (resp.\ $\alpha^{-1}(\tau)\leq 2M^{-1}$). 
\end{proof}

\begin{proof}[Proof of Theorem~\ref{thm:Z} (iii)]
	Any $Z'\in\cC(\R,\R^2)$ represents an equivalence class of processes. Let $\wt Z=(\wt Z_t)_{t\in\R}$ be an arbitrary representative for this equivalence class. Then there exists an increasing bijection $s:\R\to\R$ such that $\wt Z_t=Z_{s(t)}$ for all $t\in\R$. We want to show that $s\in\sigma(\wt Z)$, i.e., the function $s$ is measurable with respect to $\wt Z$. Since $\{t\in\R \,:\,Z_t=0\}=\{0 \}$, we know that $s^{-1}(0)\in\sigma(\wt Z)$. We may therefore assume, upon recentering $\wt Z$, that $\wt Z_0=0$, and we will make this assumption in the remainder of the proof. Since $s$ is continuous, to show that $s\in\sigma(\wt Z)$ it is sufficient to show that for any $t\in\R$, we have $\mcl L(\wt\eta([t\wedge 0,t\vee 0]))\in\sigma(\wt Z)$, where $\wt\eta(t):=\eta(s(t))$ is the reparameterized SLE curve and $\mcl L$ is Lebesgue measure. By symmetry in law of the curve under time-reversal, it is sufficient to consider the case when $t>0$. 
	
	Let $U$ be a uniform random variable with value in $[0,1]$ which is independent of $\eta$.
	For any $M\geq 1$ define stopping times $\wt T_n(M)$ for $n\in\N\cup\{0\}$ as follows
	\eqbn
	\begin{split}
		\wt T_0(M) &= \inf\{ t>0\,:\, M\wt L_t \in \{U-1,U \} \},\\
		\wt T_n(M) &= \inf\{ t>\wt T_{n-1}(M)\,:\, M|\wt L_t-\wt L_{\wt T_{n-1}(M)}|\geq 1 \}.
	\end{split}
	\eqen
	We define $\wt T_n=\wt T_n(1)$ for all $n\in\N\cup\{0 \}$. We will argue that to conclude the proof of the proposition, it is sufficient to show  that 
	\eqb
	\mcl L(\wt\eta([0,\wt T_n]))/n\rta C
	\qquad
	\text{in probability as }n\rta\infty
	\label{eq:contour10}
	\eqe
	for some deterministic constant $C$. Assume $\wt T_n$ satisfies \eqref{eq:contour10}, and choose a sequence $(\ep_k)_{k\in\N}$ converging slowly to zero, such that if $p(n)=\sup_{k\geq n}\P[|\mcl L(\wt\eta([0,\wt T_k]))/k-C|>\ep_k ]$ then $\lim_{n\rta\infty}p(n)=0$. Consider an increasing sequence $(n_k)_{k\in\N}$ such that $\sum_{k=1}^{\infty}p(n_k)<\infty$.
	For each $k\in\N$ let $M_k\in\N$ be the smallest (random) natural number such that $\wt T_{n_k}(M_k)<t$. 
	Let $\wt n_k=\sup\{ n\in\N\,:\, \wt T_n(M_k)<t \}$. By scale invariance of SLE, we have $$\P[|\mcl L(\wt\eta([0,\wt T_{\wt n_k}(M_k)]))\wt n_k^{-1}M_k^{8/5}-C|>\ep_k ]\leq p(\wt n_k),$$ and the same property holds with $\wt n_k$ replaced by $\wt n_k+1$. It follows by the Borel-Cantelli lemma that $\wt n_k C M_k^{-8/5}$ converges a.s.\ to $\mcl L(\wt\eta([0,t]))$ as $k\rta\infty$. In particular, $\mcl L(\wt\eta([0,t]))\in\sigma(Z)$. We conclude that the proposition follows once we have proved \eqref{eq:contour10}.
	
	Define the following stopping times $T_n$ for $Z$
	\eqb
	\begin{split}
		T_0 = \inf\{ t>0\,:\, |L_t-U|\not\in[0,1] \},\quad
		T_n = \inf\{ t>T_{n-1}\,:\, |L_t-L_{T_{n-1}}|\geq 1 \}.
	\end{split}
	\label{eq:contour1}
	\eqe
	Observe that $\wt T_n=s(T_n)$, which implies that $\lim_{n\rta\infty}\mcl L(\wt\eta([0,\wt T_n]))/n\rta C$ a.s.\ if and only if $\lim_{n\rta\infty} T_n/n\rta C$ a.s., so in order to complete the proof of the proposition it is sufficient to prove the latter result. By the Birkhoff ergodic theorem this will follow if we can establish the following two results, where $S_n:=T_{n}-T_{n-1}$ for any $n\in\N$: (i) $(S_n )_{n\in\N}$ is stationary and ergodic, and (ii) $\E[S_1]<\infty$.
	
	First we will establish (i). The measure on $\R$ having unit point masses at $T_n$ for $n\in\Z$ (with $T_n$ for $n<0$ defined such that the formula \eqref{eq:contour1} for $T_n$ holds also for $n\leq 0$) has a translation invariant law since $\eta$ has stationary increments (Lemma \ref{thm:ergodic}). This implies stationarity of $(S_n)_{n\in\N}$. Ergodicity of $(S_n)_{n\in\N}$ follows from Lemma \ref{lem-gfftrivial} and transience of $\eta$.
	
	Finally, we will prove (ii). It is sufficient to show that $\P[T_2>2M]$ decays super-polynomially in $M$. By a union bound,
	\eqbn
	\P[T_2>2M] \leq \P\left[ \sup_{t\in[0,M]} |L_t|<1 \right] + 
	\P\left[ \sup_{t\in[M,2M]} |L_t-L_M|<1 \right].
	\eqen
	The two terms on the right side are equal, so we will only bound the first term. For each $z\in\C$ let $\eta_z^L$ be the left frontier of $\eta$ upon hitting $z$, equipped with the natural parametrization. Defining $A_M:=\Z^2\cap[M,M]^2$ for any $M>0$ and $\tau_z:=\inf\{ t\geq 0\,:\,\eta_z^L(t)\not\in B_{M^{0.01}}(z) \}$, a union bound gives
	\eqbn
	\begin{split}
		\P\left[ \sup_{t\in[0,M]} |L_t|<1 \right]
		\leq&\,\, \P[B_{M^{0.01}}(z)\not\subset \eta([0,M]),\,\,\forall z\in A_M ]
		+ \sum_{z\in A_M} \P[\tau_z<1 ].
	\end{split}
	\eqen
	The first term on the right side decays super-polynomially in $M$ by \cite[Lemma 3.6]{kpz-bm}. The second term on the right side is decays super-polynomially in $M$ by Lemma \ref{prop10}. It follows that $\P[T_2>2M]$ decays super-polynomially in $M$, so (ii) holds.
\end{proof}

\section{The SLE is measurable with respect to the pair of contour functions}\label{sec:graph}
In this section we will prove Theorem \ref{thm:meas} and Corollary \ref{cor1}. 
For $\kappa>4$ let $\eta$ be a whole-plane space-filling SLE$_\kappa$ on $\C$ parametrized by Lebesgue measure such that $\eta(0)=0$. Recall the definition of the pair of contour functions $Z=(L,R)$ in the introduction. 
Conditioned on $Z$,  independently  sample two SLE curves $\eta$ and $\wt \eta$ according to the conditional distribution of $\eta$ given $Z$. Notice that this conditional distribution is well-defined by \cite[Theorem 5.1.9]{durrett}, since $\eta$ is a random variable with values in the space of continuous curves equipped with the topology of uniform convergence on compact sets, which is a standard Borel space. Define $\phi:\C\to\C$ by $\phi(\eta(t))=\wt\eta(t)$ for all $t\in \R$; we will prove below that $\phi$ is well-defined. In order to complete the proof of Theorem \ref{thm:meas} it is sufficient to prove the following proposition.
\begin{proposition}\label{prop:conformal}
	Almost surely $\phi$ is a rotation about the origin. 
\end{proposition}
\begin{lemma}
	For any $\kappa>4$ there is an $N\in\N$ such that the following holds a.s. The set of $t_1,t_2\in\R$ such that $\eta(t_1)=\eta(t_2)$ is exactly the set of $t_1,t_2\in\R$ for which there exists $n\leq N$ and $s_1,\dots,s_n\in\R$, such that $s_1=t_1$, $s_n=t_2$, and for any $i\in\{1,\dots,n-1\}$ we have either 
	\eqb
	R_{s_i}=R_{s_{i+1}}=\inf_{s_i\leq t\leq s_{i+1}} R_t
	\qquad\text{or}\qquad L_{s_i}=L_{s_{i+1}}=\inf_{s_i\leq t\leq s_{i+1}} L_t.
	\label{eq84}
	\eqe
	\label{prop35}
\end{lemma}
\begin{proof}
	Let $\mu$ be the measure in $\C$ associated with a $4/\sqrt{\kappa}$-quantum 
	cone independent of $\eta$, and let $\wh Z=(\wh L_t,\wh R_{t})_{t\in\R}$ 
	describe the evolution of the quantum length of the left and right, respectively, frontier of $\eta$, when we parametrize $\eta$ by quantum area. By Theorem \ref{thm:wedges} proved in \cite{wedges}, $\wh Z$ has the law of a two-dimensional correlated Brownian motion. Since $\mu$ assigns a positive measure to each open set and has no point masses, there is a continuous strictly increasing bijective function $\alpha:\R\to\R$ satisfying $\alpha(0)=0$ and $\mu(\eta([t\wedge 0,t\vee 0]))=|\alpha(t)|$ for all $t\in\R$. By the peanosphere construction of \cite{wedges} (see the text right after Corollary \ref{cor:mate}), and since for any $\kappa>4$ there is an $N\in\N$ such that no points of a space-filling SLE$_\kappa$ has multiplicity larger than $N$ (see e.g.\ \cite[Theorem 6.3]{kpz-bm}), we know that the lemma holds if we parametrize $\eta$ by quantum area and consider $\wh Z$ instead of $Z$. In other words, defining $\wh\eta(t)=\eta(\alpha^{-1}(t))$, the set of $t_1,t_2\in\R$ such that $\wh\eta(t_1)=\wh\eta(t_2)$, is exactly the set of $t_1,t_2\in\R$ for which there exists $n\leq N$ and $s_1,\dots,s_n\in\R$, such that $s_1=t_1$, $s_n=t_2$, and for any $i\in\{1,\dots,n-1\}$ one of the conditions \eqref{eq84} is satisfied. By this result and symmetry in $L$ and $R$, in order to conclude the proof of the lemma, it is sufficient to show that a.s.,
	\eqb
	\begin{split}
		&\left\{ (t_1,t_2)\in\R^2\,:\, L_{t_1}=L_{t_{2}}=\inf_{t_1\leq t\leq t_{2}} L_t \right\} \\
		&\qquad\qquad\qquad\qquad\qquad\qquad\qquad= 
		\left\{ (t_1,t_2)\in\R^2\,:\, \wh L_{\alpha(t_1)}=\wh L_{\alpha(t_{2})}=\inf_{\alpha(t_1)\leq t\leq \alpha(t_{2})} \wh L_t \right\}.
	\end{split}
	\label{eq86}
	\eqe
	Let $q\in\Q$. Almost surely, for any $t_1<q$ such that $\eta(t_1)$ is not contained in the left frontier of $\eta$ at time $q$, we have $L_{t_1}>\inf_{t\in[t_1,q]} L_t$.  Therefore, a.s., for any $t_1<q$ such that $L_{t_1}=\inf_{t\in[t_1,q]} L_t$, the left frontier of $\eta$ at time $t_1$ is contained in the left frontier of $\eta$ at time $q$. It follows that a.s.\ for any $(t_1,t_2)$ contained in the set on the left side of \eqref{eq86}, the left frontier of $\eta$ at time $t_1$ is contained in the left frontier of $\eta$ at any rational time in $[t_1,t_2]$. This implies that $\wh L_{\alpha(t_1)}=\inf_{\alpha(t_1)\leq t\leq \alpha(t_{2})} \wh L_t$. Since the time-reversal of $(R,L)$ and $(\wh R,\wh L)$ describe the evolution of the boundary lengths for the time-reversal of $\eta$, it follows that a.s.\ for any $(t_1,t_2)$ contained in the set on the left side of \eqref{eq86}, we also have $\wh L_{\alpha(t_2)}=\inf_{\alpha(t_1)\leq t\leq \alpha(t_{2})} \wh L_t$. This proves that the set on the left side of \eqref{eq86} is a.s.\ contained in the set on the right side of \eqref{eq86}. Proving that the set on the right side of \eqref{eq86} is a.s.\ contained in the set on the left side of \eqref{eq86} is done by an identical argument, and we can conclude that \eqref{eq86} holds.
\end{proof}

\begin{lemma}
	The function $\phi:\C\to\C$ defined above is a.s.\ well-defined, and $\phi$ is a.s.\ an area-preserving homeomorphism.
	\label{prop:phiwelldef}
\end{lemma}
\begin{proof}
	By Lemma \ref{prop35}, a.s.,
	\eqbn
	\{ (t_1,t_2)\in\R^2\,:\, \eta(t_1)=\eta(t_2) \} = 
	\{ (t_1,t_2)\in\R^2\,:\, \wt\eta(t_1)=\wt\eta(t_2) \}.
	\eqen
	This implies that $\phi$ is well-defined and bijective. 
	
	Next we will argue that $\phi$ is a.s.\ continuous. By symmetry in $\eta$ and $\wt\eta$, and since $\phi$ is bijective, this will imply that $\phi$ is a homeomorphism a.s. It is sufficient to argue that a.s., for any $z\in\C$, any sequence $(z_n)_{n\in\N}$ converging to $z$, and any $\ep>0$, we have $|\phi(z_n)-\phi(z)|<\ep$ for all sufficiently large $n$. Let $k\in\N$ and $t_1,\dots,t_k\in\R$ be such that $\eta(t)=z$ iff $t=t_j$ for some $j\in\{1,\dots,k\}$. For each $j\in\{1,\dots, k \}$ let $I_j$ be an open interval containing $t_j$ such that $\wt\eta(I_j)\subset B_\ep(\phi(z))$. For each $n\in\N$, let $s_n\in\R$ be such that $\eta(s_n)=z_n$. To conclude the proof, it is sufficient to argue that $s_n\in \cup_j I_j$ for all sufficiently large $n$. We will prove this by contradiction, and assume there is a strictly increasing sequence $(n_k)_{k\in\N}$ such that  $s_{n_k}\not\in \cup_j I_j$ for all $k\in\N$. The sequence $(s_n )_{n\in\N}$ is bounded, so we can find $s\in\R$ such that $(s_{n_k})_{k\in\N}$ converges subsequentially to $s$. By continuity of $\eta$ we have $\eta(s)=\lim_{k\rta\infty} \eta(s_{n_k})= \lim_{n\rta\infty} z_n=z$, so $s=t_j$ for some $j$. This is a contradiction to the assumption $s_{n_k}\not\in \cup_j I_j$, and we conclude that $\phi$ is continuous.
	
	To prove that $\phi$ is a.s.\ measure-preserving it is sufficient to prove that for any disk $B\subset\C$ we have $\mcl L(B)=\mcl L(\phi(B))$ a.s., where $\mcl L$ denotes Lebesgue measure. Let $(J_k)_{k\in\N}$ be a countable collection of disjoint intervals such that $B=\cup_{k\in\N}\eta(J_k)$. Since $\eta$ is parametrized by Lebesgue measure, $\mcl L(B)=\sum_{k\in\N}|J_k|$, where we use $|\cdot|$ to denote the length of an interval. Since $\wt\eta=\phi\circ\eta$, we know that $\phi(B)=\cup_{k\in\N}\wt\eta(J_k)$, so since $\wt\eta$ is parametrized by Lebesgue measure, we have $\mcl L(\phi(B))=\sum_{k\in\N} |J_k|$. It follows that $\mcl L(B)=\mcl L(\phi(B))$, so $\phi$ is a.s.\ measure-preserving.
\end{proof}  
For fixed $a,b\in\C$, define
\eqbn
\begin{split}
	&A_1=A_1(a,b):=|\phi(a+b)-\phi(a)|,\\ &A_2=A_2(a,b):=|\phi(a+2b)-\phi(a+b)|,\\ &A_3=A_2(a,b):=|\phi(a+2b)-\phi(a)|. 
\end{split}
\eqen
\begin{lemma}
	For any fixed $a,b\in\C$ we have
	$A_1\overset{d}{=} A_2$ and $A_3 \eqD 2 A_1$.
\end{lemma}
\begin{proof}
	We first remark that the first result $A_1\eqd A_2$ is \emph{not} an immediate consequence of invariance under recentering of whole-plane space-filling SLE, which holds by invariance under recentering of the whole-plane GFF from which the curve is generated. In order to show that $A_1\eqd A_2$, we need to show that the \emph{joint} law of $\eta$ and $\wt\eta$ is invariant under recentering the curves at the time when $\eta$ hits $b$.
	
	Let $X = (Z,\eta,\wt\eta)$ be the triple consisting of the boundary length process $Z$, in addition to the two conditionally independent curves $\eta$ and $\wt\eta$. For any $t\in\R$, let $X(t)$ be equal to $X$, except that the processes are recentered at time $t$, i.e., 
	\eqbn
	X(t) = (Z_{\cdot+t}-Z_t,
	\eta(\cdot+t)-\eta(t), 
	\wt\eta(\cdot+t) -\wt\eta(t)
	).
	\eqen
	Fix $M>0$, and let $\sigma$ be a uniform random variable in $[-M,M]$ independent of $X$. For any $z\in\C$, let $\tau(z)=\inf\{t\in\R\,:\,\eta(t)=z \}$ be the time at which $\eta$ hits $z$. By independence of $\sigma$ and $X$, and since $X\eqD X(t)$ for any fixed $t\in\R$ by Lemma \ref{thm:ergodic}, we have $X\eqD X(\sigma)$, so
	\eqb
	(X,X(\tau(b))) \eqD 
	(X(\sigma) , X(\tau(\eta(\sigma)+b))).
	\label{eq81}
	\eqe
	When $M\rta\infty$, the law of $X(\tau(\eta(\sigma)+b))$ converges to the law of $X(\sigma)=X(\tau(\eta(\sigma)))$ in total variation distance, since sampling a time $\sigma$ uniformly from $[-M,M]$ is equivalent to sampling a point $z$ uniformly from $\eta([-M,M])$ in the sense that $z\eqD\eta(\sigma)$. Since the total variation distance between the laws of the two elements on the right side of \eqref{eq81} converges to zero when $M\rta\infty$, we see that the two elements on the left side of \eqref{eq81} are equal in law. This implies that $A_1\eqD A_2$.
	
	Next we will prove that $A_3 \eqD 2 A_1$. Since $X\eqD X(\tau(a))$ by the arguments of the preceding paragraph, we may assume $a=0$ in the remainder of the proof. Define $\eta'=2\eta(\cdot/4)$, 
	$\wt \eta'=2\wt \eta(\cdot/4)$, and $Z'_t=2^{ 1+2/\kappa}Z_{\cdot/4}$. Define $\acute{\phi}:\C\to\C$ such that $\acute{\phi}(\eta'(t))=\wt{\eta}'(t)$ for all $t\in\R$. Since $(Z',\eta',\wt{\eta}')\eqD (Z,\eta,\wt{\eta})$, we see that $\acute{\phi}$ is well-defined and $\acute{\phi}\overset{d}{=}\phi$. Then define $A_3'=|\acute{\phi}(2b)-\acute{\phi}(0)|$. Since $(\eta',\wt{\eta}')\eqD (\eta,\wt{\eta})$
	and $A_3'=2 A_1$, we have $A_3\eqD A_3'=2A_1$, and the second identity follows.
\end{proof}

The following proposition will be proved in Section~\ref{sec:expectation}.
\begin{proposition}\label{prop:expectation}
	For fixed $a,b\in\C$, we have $\E(A_1)<\infty.$
\end{proposition}
Combining the above results we can show that $\phi$ is linear.
\begin{proposition}	\label{prop:philinear}
	The map $\phi:\C\to\C$ is a.s.\ linear, and the matrix describing the linear transformation has determinant $\pm 1$.
\end{proposition}
\begin{proof}
	Let $a,b\in\C$. Since $\E(A_1(a,b))<\infty$ and $A_3(a,b) \overset{d}{=} 2 A_1(a,b)\eqD 2A_2(a,b)$, we have $\E[A_1(a,b)]+\E[A_2(a,b)]=\E[A_3(a,b)]$. Since we know by the triangle inequality that $A_1(a,b)+A_2(a,b)\geq A_3(a,b)$, we see that $A_1(a,b)+A_2(a,b)=A_3(a,b)$ a.s.
	This means that $\phi(a),\phi(a+b),\phi(a+2b)$ are a.s.\ collinear. Therefore, a.s.\ for any $q,\wt q\in\Q$ and $m,k\in\Z$, the following three points are collinear
	\eqb
	\phi(\wt q+ q i ),\quad
	\phi(\wt q+ q i+m 2^{k} ),\quad
	\phi(\wt q+ q i+m 2^{k+1} ).
	\label{eq87}
	\eqe
	Furthermore, the second point of \eqref{eq87} is a.s.\ between the first point and the third point of \eqref{eq87}.
	
	For $q\in\R$ define $\ell_q:=\{ z\in\C\,:\,\op{Re}(z)=q \}$ 
	and
	$\wt\ell_q:=\{ z\in\C\,:\,\op{Im}(z)=q \}$.
	By continuity of $\phi$, and since the three points \eqref{eq87} are collinear, for any fixed $q$ the set $\phi(\ell_q)$ is a.s.\ contained in a line. By Lemma \ref{prop:phiwelldef}, $\phi(\ell_q)$ is homeomorphic to $\ell_q$, so $\phi(\ell_q)$ is a.s.\ either a line segment, a half-line, or a line. Furthermore, by continuity of $\phi$ we know that this property holds a.s.\ simultaneously for all lines $\ell$, since a.s., any three collinear points are mapped to three collinear points. By symmetry in $\eta$ and $\wt\eta$, we know that $\phi^{-1}$ also maps any line to a line segment, a half-line, or a line. Using that $\phi$ is bijective, this implies that $\phi$ (and $\phi^{-1}$) maps any line to a line.
	
	For any given $k\in\Z$ consider the grid made by the lines $\ell_{m 2^k}$ and $\wt\ell_{m 2^k}$ for $m\in\Z$. Since $\phi$ is a homeomorphism a.s., any connected domain $D$ bounded by four of the grid lines $\ell_{m 2^k}$ and $\wt\ell_{m 2^k}$, is mapped bijectively onto the domain bounded by the image of these lines under $\phi$. Since $\phi$ is bijective, we see that the lines $\phi(\ell_{m 2^k})$ for different $m$ cannot intersect, so they are parallel, and the same property holds for the lines $\wt \ell_{m 2^k}$. We conclude that a.s.\ for any fixed $k\in\Z$, $\phi$ is an affine map restricted to the grid made by $\ell_{m 2^k}$ and $\wt\ell_{m 2^k}$ for $m\in\Z$, so since $\phi(0)=0$, $\phi$ is continuous, and $k$ was arbitrary, the map $\phi$ must be linear. The determinant of the matrix describing the linear map must be $\pm 1$ since the map is measure-preserving.
\end{proof}

\begin{figure}[ht!]
	\begin{center}
		\includegraphics[scale=0.75]{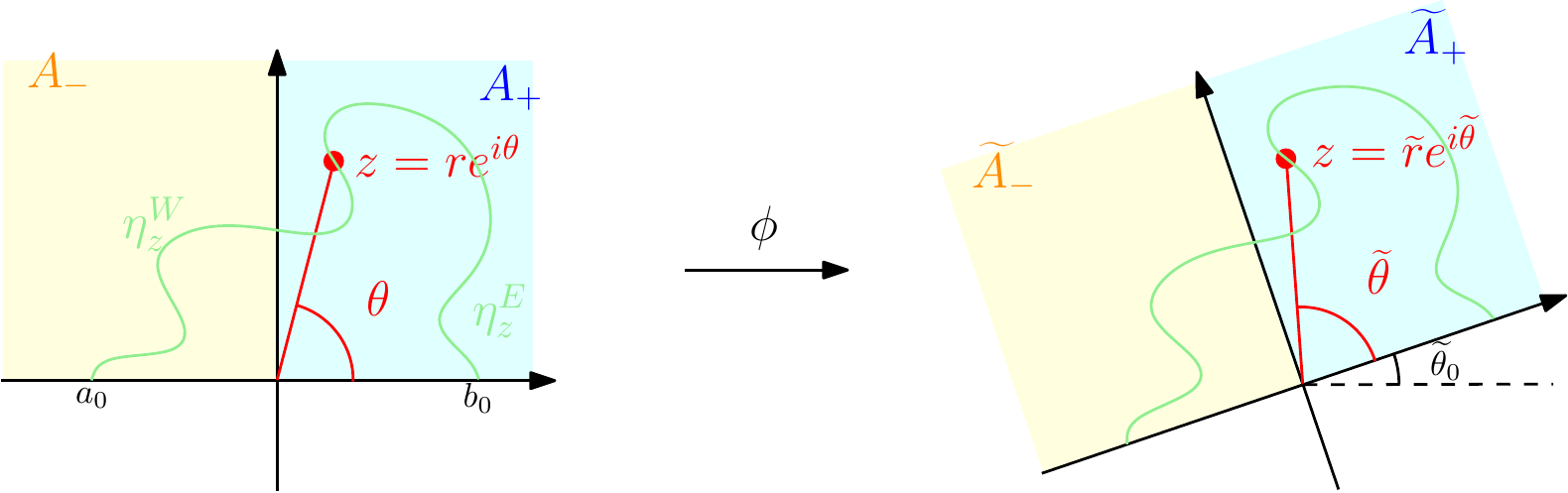} 
		\caption{
			Illustration of the proof of Proposition \ref{prop:conformal}. The region $A_+$ (resp.\ $A_-$) in light blue (resp.\ yellow) on the left figure corresponds to the points $(\theta_0,\theta,r)$ with $\theta_0=0$ for which $\lim_{n\rta\infty} Z^{\theta_0,\theta,r}_n>\frac12$ (resp.\ $\lim_{n\rta\infty} Z^{\theta_0,\theta,r}_n<\frac12$). The right figure shows $\wt A_\pm$, which are a.s.\ identical to $\phi( A_\pm)$ by the definition of $\phi$.   Since $\phi\circ\eta$ has the marginal law of an SLE, each of the domains $\wt A_\pm$ is a rotation (possibly composed with a reflection) of the first quadrant a.s.   Since a reflection would interchange the two coordinates of $Z$, we conclude that $\phi$ is a rotation a.s.
		} \label{fig-coordinates}
	\end{center}
\end{figure}

The following lemma is the main ingredient used to deduce Proposition \ref{prop:conformal} from Proposition \ref{prop:philinear}. See Figure \ref{fig-coordinates} for an illustration. Using that both $\eta$ and $\wt\eta$ have the marginal law of an SLE$_\kappa$, we will use the lemma to deduce that the linear map $\phi$ preserves angles. It will be proved in Section~\ref{sec:expectation}.
\begin{lemma}\label{lem:ergodicity}
	Let $\ul\kappa<4$, let $h$ be a whole-plane GFF modulo $2\pi\chi$ with $\chi$ given by \eqref{eq92}, and for each $z\in\C$ let $\eta_z^L$ (resp.\ $\eta_z^R$) be the west-going (resp.\ east-going) flow line of $h$ started at $z$ with some arbitrary  monotone parametrization. For $\theta\in(0,\pi)$, $r>0$, and $k\in \N$ define $z_k:=2^{-k}re^{i\theta}$,
	\eqbn
	\begin{split}
		&t_k^L := \inf\{ t\geq 0 \,:\, \op{Im} \eta^L_{z_k}(t)<0 \},\,\,\,\,
		t_k^R := \inf\{t\geq 0\,:\, \op{Im}\eta^R_{z_k}(t)<0 \},\\
		&a_k:= \eta^L_{z_k}(t_k^L),\,\,\,\,
		b_k:= \eta^R_{z_k}(t_k^R),\,\,\,\,
		x_k:=\1_{\{a_k+b_k\ge 0\}},\,\,\,\,
		Z_n:=\frac 1n \sum_{k=0}^{n}x_k.
	\end{split}
	\eqen
	For $\theta<\pi/2$ (resp.\ $\theta=\pi/2$, $\theta>\pi/2$) there is a constant $p_\theta\in(0,1)$ satisfying $p_\theta>0.5$ (resp. $p_\theta=0.5$, $p_\theta<0.5$), such that a.s.- $\lim_{n\to\infty}Z_n=p_\theta$.
\end{lemma}

We conclude the proof of Proposition~\ref{prop:conformal} by showing that unless $\phi$ is of the desired form, $\eta$ and $\wt\eta$ cannot both satisfy the property of Lemma \ref{lem:ergodicity}.
\begin{proof}[Proof of Proposition~\ref{prop:conformal}]
	We will argue that the image of a pair of orthogonal lines is a.s.\ mapped to a pair of orthogonal lines under $\phi$. This is sufficient to complete the proof, since it implies by Proposition \ref{prop:philinear} that $\phi$ is a composition of a rotation and possibly a reflection, and we see that $\phi$ has to be a rotation (not composed with a reflection), since the boundary length process $Z$ is invariant under a rotation of $\eta$, while the two coordinates are swapped upon a reflection.
	
	For any $\theta_0\in[0,2\pi)$, $\theta\in(0,\pi)$ and $r>0$, let $Z_n^{\theta_0,\theta,r}$ denote the random variable $Z_n$ defined in Lemma \ref{lem:ergodicity} for the curve $(e^{-i\theta_0}\eta(t))_{t\in\R}$. Define  
	\eqb
	\begin{split}
		A_-&:=\{(\theta_0,\theta,r)\in[0,2\pi)\times(0,\pi)\times\R_+\,:\,\lim_{n\rta\infty} Z_n^{\theta_0,\theta,r} < 1/2\},\\
		A_+&:=\{(\theta_0,\theta,r)\in[0,2\pi)\times(0,\pi)\times\R_+\,:\,\lim_{n\rta\infty} Z_n^{\theta_0,\theta,r} > 1/2\}.
	\end{split}
	\label{eq83}
	\eqe
	By Fubini's theorem, Lemma \ref{lem:ergodicity}, and rotational invariance of whole-plane space-filling SLE, we have $\mcl L(A_\pm\Delta B_\pm)=0$ a.s.\ for $B_+:=[0,2\pi)\times(0,\pi/2)\times \R_+$ and $B_-:=[0,2\pi)\times(\pi/2,\pi)\times \R_+$, where $\mcl L$ denotes Lebesgue measure and $\Delta$ denotes symmetric difference. Recall that $\wt\eta=\phi\circ\eta$, and let $\wt A_\pm$ be defined exactly as $A_\pm$, but with $\wt\eta$ instead of $\eta$, i.e., we first define $\wt Z^{\theta_0,\theta,r}_n$ exactly as $Z^{\theta_0,\theta,r}_n$ with $\wt\eta$ instead of $\eta$, and then we define $\wt A_\pm$ by \eqref{eq83} with $\wt Z^{\theta_0,\theta,r}_n$ instead of $Z^{\theta_0,\theta,r}_n$. By definition of $\phi$, and since for any $a,b\in\C$ for which $a,b,0$ are collinear, we have $|a|>|b|$ iff $|\phi(a)|>|\phi(b)|$, we have 
	\eqbn
	\begin{split}
		\wt A_\pm=& \{(\wt\theta_0,\wt\theta,\wt r)\in [0,2\pi)\times(0,\pi)\times\R_+\,:\,\exists (\theta_0,\theta,r) \in A_\pm \text{\, such\,\,that\, } \\
		&\wt\theta_0=\op{arg} \phi(e^{i\theta_0}),\,\, 
		\wt re^{i(\wt\theta_0+\wt\theta)}=\phi(re^{i(\theta_0+\theta)}) \}.
	\end{split}
	\eqen 
	Since $\wt\eta$ has the marginal law of an SLE$_\kappa$, we see from Lemma \ref{lem:ergodicity} that 
	\eqb
	\mcl L(\wt A_\pm\Delta B_\pm)=0.
	\label{eq10}
	\eqe
	In the remainder of the proof we assume that the matrix describing the linear transformation $\phi$ has determinant 1 (equivalently, the curve $(\phi(e^{i\theta}))_{\theta\in[0,2\pi)}$ goes counterclockwise about the origin); the opposite case can be treated similarly.
	Let $\theta_0$ be sampled uniformly at random from $[0,2\pi)$. Letting
	\eqbn
	\begin{split}
			&A_\pm^{\theta_0} = \{ (\theta,r)\in(0,\pi)\times\R_+\,:\, (\theta_0,\theta,r)\in A_\pm \},\\
		& 
		B^{\pi/2}_+=(0,\pi/2)\times\R_+,
		\qquad
		B^{\pi/2}_-=(\pi/2,\pi)\times\R_+,
	\end{split}
	\eqen
	it follows from $\cL(A_\pm\Delta B_\pm)=0$ that a.s.,  $\cL(A^{\theta_0}_\pm\Delta B^{\pi/2}_\pm)=0$. Let $\theta^*_0\in[0,2\pi)$ be the angle between the positive $x$-axis and the image of $s=\{ re^{i\theta_0}\,:\, r\geq 0 \}$ under $\phi$ in counterclockwise direction, and let $\theta$ be the angle between the images of $s=\{ re^{i\theta_0}\,:\, r\geq 0 \}$ and 
	$s'=\{ 
	re^{i(\theta_0+\pi/2)}\,:\, r\geq 0
	\}$ under $\phi$ in counterclockwise direction. We have $\theta^*\in(0,\pi)$ by our assumption that the determinant of the matrix describing $\phi$ is equal to 1. Defining
	\eqbn
\begin{split}
	&\wt A_\pm^{\theta^*_0} = \{ (\theta,r)\in(0,\pi)\times\R_+\,:\, (\theta^*_0,\theta,r)\in \wt A_\pm \},\\
	& 
	\wt B^{\theta^*}_+=(0,\theta^*)\times\R_+,
	\qquad
	\wt B^{\theta^*}_-=(\theta^*,\pi)\times\R_+,
\end{split}
\eqen	
it follows from $\cL(A^{\theta_0}_\pm\Delta B^{\pi/2}_\pm)=0$ that $\cL(\wt A^{\theta^*_0}_\pm\Delta \wt B^{\theta^*}_\pm)=0$ a.s. On the other hand, it follows from \eqref{eq10} that $\cL(\wt A^{\theta^*_0}_\pm\Delta B^{\pi/2}_\pm)=0$. This implies that $B^{\pi/2}_\pm=\wt B^{\theta^*}_\pm$ a.s., so $\theta^*=\pi/2$ a.s.\ and the two orthogonal lines $s$ and $s'$ are mapped to orthogonal lines under $\phi$. 
\end{proof}

\begin{proof}[Proof of Corollary \ref{cor1}]
	It is sufficient to prove that we get a topological sphere when we glue together the pair of trees in Figure \ref{fig-mating}; once we have proved this it is immediate that the sphere is equipped with a space-filling path (mapping each $t_0\in\R$ to the equivalence class of $(\phi(t_0),R'_{\phi(t_0)})$), and Theorem \ref{thm:meas} implies that the embedding of the path-decorated sphere into $\C$ is canonical. It follows by Lemma \ref{prop35} that we get a topological sphere under the equivalence relation on Figure \ref{fig-mating}, since the existence of appropriate times $s_1,\dots,s_k$ is exactly the condition which says whether two times $t_1$ and $t_2$ are in the same equivalence class for the considered equivalence relation, and since the lemma implies that the quotient topology on the set of equivalence classes is the same as the standard topology on $\C$.
\end{proof}

\section{Proof of Proposition~\ref{prop:expectation} and Lemma~\ref{lem:ergodicity}}
\label{sec:expectation}
The proofs of Proposition~\ref{prop:expectation} and Lemma~\ref{lem:ergodicity} are based on regularity estimates for space-filling SLE$_\kappa$, which we will prove in Lemma \ref{lem:exp1} and Proposition \ref{prop:powerlaw}. Throughout the section let $\kappa>4$, let $\chi$ be given by \eqref{eq92} with $\ul\kappa=16/\kappa$, let $h$ be a whole-plane GFF modulo a global additive multiple of $2\pi\chi$, and for each $z\in\C$ let $\eta^L_z$ (resp.\ $\eta^R_z$) be the west-going (resp.\ east-going) flow line of $h$ started from $z$. Let $\eta$ be the whole-plane space-filling SLE$_\kappa$ generated by $h$, parametrized by Lebesgue measure and satisfying $\eta(0)=0$.

\begin{figure}[ht!]
	\begin{center}
		\includegraphics[scale=0.8]{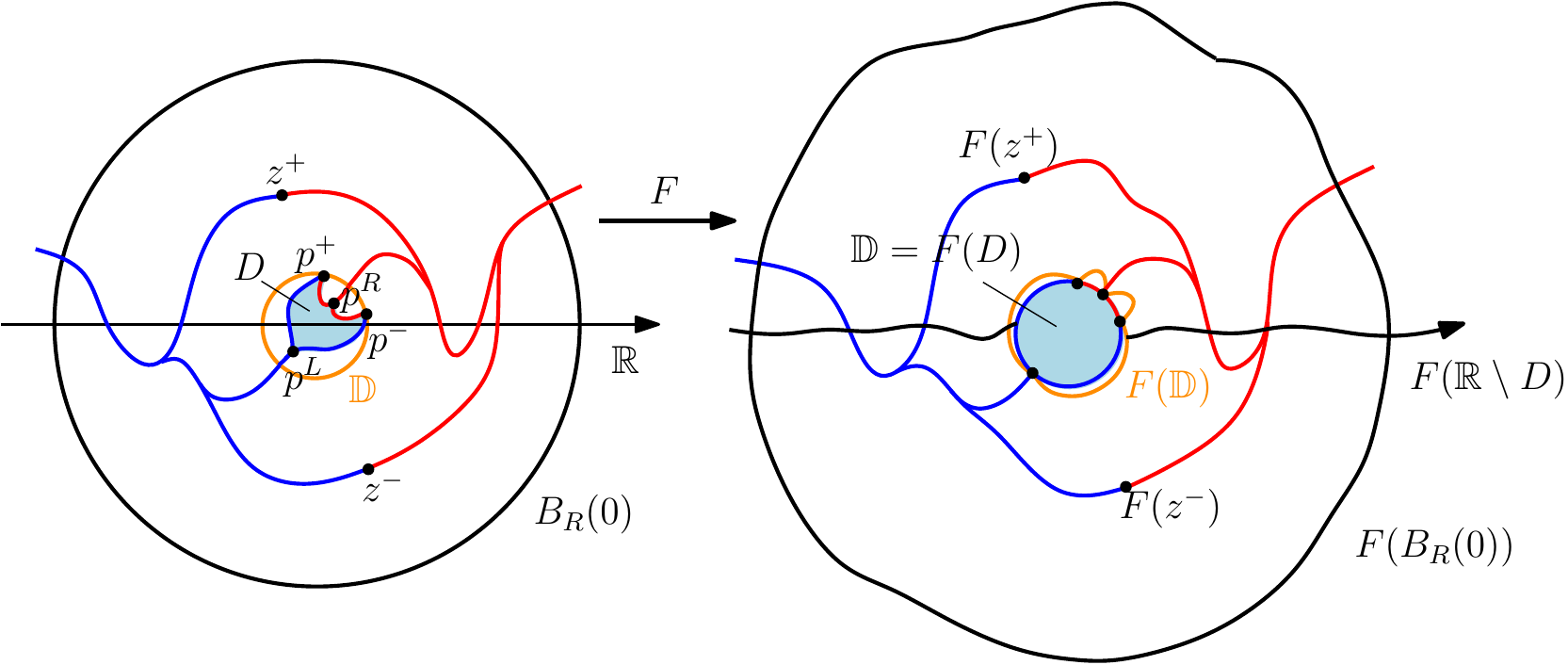} 
		\caption{The left figure illustrates the event $E_R$ in the statement of Lemma \ref{lem:exp1}, and the right figure illustrates the proof of the lemma. We show that $E_R$ has a uniformly positive probability of occuring conditioned on the set $D$ (shown in light blue) and $\frk p=(p^+,p^-,p^R,p^L)$. We do this by applying Lemma \ref{lem:RN}, which says  that the realization of the Gaussian free field on $\C\setminus\D$ on the right figure does not depend too strongly on $D$ and $\frk p$ in domains bounded away from $\D$ and $\infty$.}
		\label{fig:pictureI}
	\end{center}
\end{figure}

The lemma we state next will be applied in the proof of both Proposition~\ref{prop:powerlaw} and Lemma~\ref{lem:ergodicity}. Define stopping times $\sigma^\pm$ for $\eta$ as follows 
\begin{equation}
\sigma^+=\sup\{t\geq 0:\eta([0,t])\subset \D \},\quad
\sigma^-=\inf\{t\leq 0:\eta([t,0]) \subset \D\}.
\end{equation}
Let $p^\pm:=\eta(\sigma^\pm)$. Let $p^R$ (resp.\ $p^L$) be the point at which the two flow lines $\eta_{p^\pm}^R$ (resp.\ $\eta_{p^\pm}^L$) merge. Then define the $\sigma$-algebra $\mcl G$ by $\mcl G:=\sigma(D,\frk p)$, where $\frk p\in\C^4$ and $D\subset\C$ are defined by
\eqbn
\begin{split}
	&\frk p:=(p^+,p^-,p^R,p^L),\quad 
	D:=\eta([\sigma^-,\sigma^+]).
\end{split}
\eqen
See the left part of Figure \ref{fig:pictureI} for an illustration of the objects defined above and of the statement of the following lemma.
\begin{lemma}
	In the setting described above, for $R>1$ let $E_R$ be the event that there exists $z^+,z^-\in B_R(0)$ such that $z^+$ (resp.\ $z^-$) is contained in the upper (resp.\ lower) half-plane, and such that the following hold
	\begin{itemize}
		\item $\eta^L_{z^+}$ and $\eta^L_{z^-}$ (resp.\ $\eta^R_{z^+}$ and $\eta^R_{z^-}$) merge before they exit $B_R(0)$,
		\item $\eta^L_{z^+}$ and $\eta^L_{z^-}$ (resp.\ $\eta^R_{z^+}$ and $\eta^R_{z^-}$) hit $\R_-$ (resp.\ $\R_+$) before they exit $B_R(0)$, and
		\item the bounded region enclosed by the four flow lines $\eta^L_{z^+},\eta^L_{z^-},\eta^R_{z^+},\eta^R_{z^-}$ contains $\D$.
	\end{itemize}
	Then there exist $p>0$ and $R_0>0$, such that for $R\geq R_0$ we have $\P[E_R\,|\,\mcl G]> p$.
	\label{lem:exp1}
\end{lemma}

\begin{proof}
	Since the event $E_R$ is monotone in $R$, it is sufficient to prove that there exist $R>0$ and $p>0$ such that $\P[E_R\,|\,\mcl G]\ge p$.
	Let $F:\C\setminus D\to \C\setminus\D$ be the unique conformal map such that $\lim_{z\rta\infty} F(z)/z>0$. The logarithmic capacity of $D$, which is denoted by cap$(D)$, is defined to be 
cap$(D):=\log\lim_{z\rta\infty} F(z)/z$; see e.g. \cite[Chapter 3]{lawler-book}. By \cite[Propositions 3.29-3.30]{lawler-book} there is a universal constant $c>0$ such that $|F(z)-e^{-\op{cap}(D)}z|<c$ for all $z\in\C\setminus D$, where $\op{cap}(D)$ is the logarithmic capacity of $D$ and $e^{-\op{cap}(D)}\in[1,4]$. Let $R=100(1+c)$, and observe that 
	\eqb
	\begin{split}
		&F(\D)\subset B_{0.1 R}(0),\qquad B_{0.9R}(0)\subset F(B_R(0))\\
		&F(\R\setminus D) \subset \{z\in\C\,:\,|\op{Im}(z)|<0.1R \}.
	\end{split}
	\label{eq:exp1}
	\eqe
	
	We will argue that $D$ is a local set for the GFF $h$, as defined in \cite[Section 3.3]{ss-contour} (see also \cite[Section 3.2]{IG1}). Given any $z\in\C$ and $\theta\in[0,2\pi)$ the flow line $\eta_z^\theta$ of $h$ started at $z$ with angle $\theta$ is a local set for $h$ by \cite[Theorem 1.1]{IG4} and \cite[Lemma 3.9, 4.]{ss-contour}. Let $H(\C)$ be the Hilbert space closure for the Dirichlet inner product of the space $C^\infty_c(\C)$ of real-valued smooth compactly supported functions on $\C$. For any open set $U\subset\C$, let $H_{\op{supp}}(U)$ be the subset of $H(\C)$ consisting of functions which are supported in $U$, and let $H_{\op{harm}}(U)\subset H(\C)$ be the orthogonal complement of $H_{\op{supp}}(U)$ for the Dirichlet inner product. Since $D$ is measurable with respect to a countable collection of flow lines for $h$, the event that $U\cap D=\emptyset$ is measurable with respect to the projection of $h$ onto $H_{\op{harm}}(U)$, so $D$ is local by \cite[Lemma 3.9, 1.]{ss-contour}.
	
	We will now describe the boundary conditions of $h|_{\C\setminus D}$. By \cite{IG4}, for any fixed $z\in\C$ and $\theta\in[0,2\pi)$ and with $\chi$ and $\lambda'$ as in \eqref{eq92}, the boundary conditions modulo $2\pi\chi$ on the left (resp.\ right) side of $\eta_z^\theta$ are given by $-\lambda'-\theta\chi$ (resp.\ $\lambda'-\theta\chi$), plus $\chi$ times the winding of the curve, where the winding is defined relative to a path going straight upwards (equivalently, straight northwards, or in the direction of the positive imaginary axis). We say that the flow line has flow line boundary conditions $-\lambda'-\theta\chi$ (resp.\ $\lambda'-\theta\chi$) on its left (resp.\ right) side. See \cite[Theorem 1.1 and Figure 1.9]{IG4}. The flow line boundary conditions of $h|_{\C\setminus D}$ are therefore given by $-\lambda'-\frac{\pi}{2}\chi$ (resp.\ $\lambda'-\frac{\pi}{2}\chi$, $-\lambda'+\frac{\pi}{2}\chi$, $+\lambda'+\frac{\pi}{2}\chi$) on the segment $p^-p^L$ (resp.\ $p^+p^L$, $p^+p^R$, $p^-p^R$) of $\partial D$. Since $D$ is local, the characterization of local sets in \cite[Lemma 3.9]{ss-contour} implies that the conditional law of $h|_{\C\setminus D}$ given $\mcl G$ is that of a zero boundary GFF plus the harmonic extension of the values of $h_{\C\setminus D}$ from $\partial D$ to $\C\setminus D$.
	
	Define the distribution $\wh h$ on $\C\setminus \D$ by $\wh h:= h\circ F^{-1}-\chi\op{arg}(F^{-1})'$. For any $z\in\C$ and $\theta\in[0,2\pi)$, the curve $F\circ \eta_z^\theta$ is a flow line for $\wh h$ of angle $\theta$ started at $F(z)$, see \cite[Section 1.3]{IG1}. Define $\wh p^\bullet:=F(p^\bullet)$ for $\bullet\in\{+,-,e,w \}$. Observe that $\wh h$ has the same flow line boundary values as $h$, i.e., the flow line boundary conditions are given by $-\lambda'-\frac{\pi}{2}\chi$ on the segment $\wh p^-\wh p^L$ of $\partial\D$, etc.
	
	Let $\wh E_R$ be the event that we can find $\wh z^+,\wh z^-\in B_{0.9R}(0)\setminus\ B_{0.1R}(0)$ satisfying $\op{Im}(\wh z^+)>0.1R$ and $\op{Im}(\wh z^-)<-0.1R$ such that
	\begin{itemize}
		\item the flow lines $\eta^L_{\wh z^+}$ and $\eta^L_{\wh z^-}$ (resp.\ $\eta^R_{\wh z^+}$ and $\eta^R_{\wh z^-}$) for $\wh h$ merge before they exit $B_{0.9R}(0)$,
		\item the flow lines $\eta^L_{\wh z^\pm}$ (resp.\ $\eta^R_{\wh z^\pm}$) hit $\R_-\mp 0.1R$ (resp.\ $\R_+\mp 0.1R$) before they exit $B_{0.9R}(0)$, and
		\item the bounded region enclosed by the four flow lines $\eta^L_{\wh z^+},\eta^L_{\wh z^-},\eta^R_{\wh z^+},\eta^R_{\wh z^-}$ for $\wh h$ contains $B_{0.1R}(0)$.
	\end{itemize}
	Observe that $\wh h|_{\partial \D}$ is bounded, and that $\wh E_R$ is measurable with respect to $\wh h|_{B_{0.9R}(0)\setminus\ B_{0.1R}(0)}$. Since the event $\wh E_R$ occurs with positive probability for any fixed choice of boundary data for $\wh h$ on $\partial\D$, and since the boundary data of $\wh h$ are bounded, it follows by Lemma \ref{lem:RN} that there is a $p>0$ such that $\P[\wh E_R\,|\,\mcl G]> p$.
	By \eqref{eq:exp1} and our choice of $R$, we have $\wh E_R\subset E_R$, since we can define $z^\pm=F^{-1}(\wh z^\pm)$ on the event that $\wh E_R$ occurs. It follows that $\P[E_R\,|\,\mcl G]>p$.
\end{proof}

\begin{proposition}\label{prop:powerlaw}
	Let $\eta$ be a whole-plane space-filling $\SLE_{\kappa}$ for $\kappa>4$, parametrized by Lebesgue measure and satisfying $\eta(0)=0$. Then there exist a $\xi>0$ such that for all $M>0$,
	\begin{equation}
	\P[\D\nsubseteq \eta([-M,M])] \preceq M^{-\xi},
	\end{equation}
	where the implicit constant may depend on $\kappa$.
\end{proposition}
\begin{proof}
	Let $E_R$ be the event of Lemma \ref{lem:exp1}, and fix $R>0$ sufficiently large such that $\P[E_R\,|\,\mcl G]>p$ for some $p>0$. For $k\in\N$ let $E_R^k$ be the event that $E_R$ holds for the Gaussian free field $h\circ g_k-\chi \op{arg}(g'_k)=h\circ g_k$, where $g_k(z):=R^{k-1}z$. In other words, $E_R^k$ is defined exactly as $E_R$, except that $B_{R}(0)$ is replaced by $B_{R^{k}}(0)$ and $\D$ is replaced by $B_{R^{k-1}}(0)$. Let $\mcl G_k$ be the $\sigma$-algebra which is defined exactly as $\mcl G$, but for the Gaussian free field $h\circ g_k$. By conformal invariance of $h$ and Lemma \ref{lem:exp1}, $\P[E_R^k\,|\,\mcl G_k]> p$ for all $k\in\N$, so $\P[\cap_{1\leq k\leq K} (E_R^k)^c]<(1-p)^K$.
	Observe that if $M=\pi R^{2K}$ for some $K\in\N$ then $\{\D \not\subset \eta([-M,M]) \}\subset \cap_{1\leq k\leq K} (E_R^k)^c$, so
	\eqbn
	\P[\D\not\subset \eta([-M,M]) ] < (1-p)^K,\quad M=\pi R^{2K}.
	\eqen
	This implies the existence of an appropriate $\xi$.
\end{proof}

\begin{proof}[Proof of Proposition~\ref{prop:expectation}]
	We will show that $\P[|\phi(a+b)-\phi(a)|>k]$ decays faster than any negative power of $k$, which is sufficient to complete the proof of the proposition. When proving this, we will consider an infinite graph $G$ defined as follows. Each vertex of $G$ is identified with an interval of the form $[m,m+1]$ for $m\in\Z$. There is an edge between vertices corresponding to intervals $[m_1,m_1+1]$ and $[m_2,m_2+1]$ iff $\eta([m_1,m_1+1])\cap \eta([m_1,m_1+1])\neq\emptyset$. We remark that $G$ is defined similarly as the \emph{structure graphs} considered in \cite{ghs-metric}, where the graphs were used to define a discrete metric on a Liouville quantum gravity surface. We note that $\eta$ and $\wt\eta$ give the same graph $G$, since $G$ is measurable with respect to $Z$ by Lemma \ref{prop35}.
	
	We fix $K>0$, and want to show that $\P[|\phi(a+b)-\phi(a)|>k]\preceq k^{-K}$ for all $k\geq 10^{10}(1+b^{10})$, where the implicit constant may depend on $K$, but not on $k$. Let $d\in\N$ be the number of vertices $[m,m+1]$ of $G$ for which $\eta([m,m+1])\cap[a,a+b]\neq\emptyset$, where $[a,a+b]$ denotes the line segment connecting $a$ and $a+b$. For any $K'>0$, a union bound gives
	\eqb
	\begin{split}
		\P\big(|\phi(a+b)-&\phi(a)|>k\big) \leq\,\, 
		\P\big([a,a+b]\not\subset \eta([-k^{K'},k^{K'}])\big)\\
		&+ \P\big([a,a+b]\subset \eta([-k^{K'},k^{K'}]);\,d\geq k^{1/2}
		\big)\\
		&+ \P\big([a,a+b]\subset \eta([-k^{K'},k^{K'}]);\,d<k^{1/2};\,
		|\phi(a+b)-\phi(a)|>k\big).
	\end{split}
	\label{eq88}
	\eqe
	Choose $K'$ sufficiently large such that the first term on the right side of \eqref{eq88} is $\preceq k^{-K}$; such a value of $K'$ exists by Proposition~\ref{prop:powerlaw}. If $d\geq k^{1/2}$, there are $\geq k^{1/2}$ cells of area 1 which intersect $[a,a+b]$, hence at least one of the cells has diameter larger than $k^{1/10}$; otherwise all the $\ge k^{1/2}$ cells would be contained in the ball $B_{k^{1/10}+b+1}(a)$, thus contradicting the fact that  the area of $B_{k^{1/10}+b+1}(a)$ is smaller than $k^{1/2}$. By a union bound, translation invariance in law of $\eta$, and \cite[Lemma 3.6]{kpz-bm},
	\eqbn
	\begin{split}
		\P\big([a,a+b]\subset \eta([-k^{K'},k^{K'}]);\,d\geq k^{1/2}\big)
		&\leq \sum_{j=-k^{K'}}^{k^{K'}-1} \P\big( \mathrm{diam}(\eta([j,j+1]))> k^{1/10}\big)\\
		&\leq 2k^{K'}\P\big( \mathrm{diam}(\eta([0,1]))> k^{1/10}\big)\preceq k^{-K}.
	\end{split}
	\eqen
	If the event in the third term on the right side of \eqref{eq88} occurs, there is an $m\in\{ -k^{K'},-k^{K'}+1,\dots, k^{K'}-1 \}$ such that $\op{diam}(\wt\eta([m,m+1]))\geq |\phi(a+b)-\phi(a)|/d>k^{1/2}$. Applying \cite[Lemma 3.6]{kpz-bm} again, we get
	\eqbn
	\begin{split}
		\P\big([a,a+b]\subset\eta([-k^{K'},k^{K'}]);\,d<k^{1/2};\,
		&|\phi(a+b)-\phi(a)|>k\big)\\
		&\leq 2k^{K'}\P\big( \mathrm{diam}(\wt\eta([0,1]))> k^{1/2}\big)\preceq k^{-K}.
	\end{split}
	\eqen
	Combining the above bounds, we see from \eqref{eq88} that $\P[|\phi(a+b)-\phi(a)|>k]\preceq k^{-K}$, which concludes the proof of Proposition~\ref{prop:expectation}.
\end{proof}

\begin{proof}[Proof of Lemma \ref{lem:ergodicity}]
	By scale invariance of SLE, it is sufficient to prove the lemma for $r=1$, and we will make this assumption throughout the proof. Define $p_\theta:=\P(x_1=1)$. We will argue that $p_\theta$ satisfies the inequality in the statement of the lemma. By symmetry, we have $p_{\pi/2}=0.5$, and we assume for the remainder of the paragraph that $\theta\neq\pi/2$. By invariance under recentering of the whole plane GFF from which the flow lines $\eta^R_{z_1}$ and $\eta^L_{z_1}$ are generated, the law of $a_k+b_k$ is symmetric about $\op{Re}(z_1)$. It holds with positive probability that $a_k+b_k\in(\op{Re}(z_1)\wedge 0,\op{Re}(z_1)\vee 0)$. Since $\theta<\pi/2$ (resp.\ $\theta>\pi/2$) implies that $\op{Re}(z_1)>0$ (resp.\ $\op{Re}(z_1)<0$) it follows that $p_\theta$ satisfies the inequalities in the statement of the lemma.
	
	First we will prove that we can find a $c>0$ such that
	\eqb
	\P(E_k^c)\leq \exp(-ck),\qquad
	E_k:=\left\{ \eta^R_{z_k}([0,t_k^R]) \subset \D;\,\eta^L_{z_k}([0,t_k^L]) \subset \D \right\}.
	\label{eq:conf1}
	\eqe
	For $k\in\N$ let $\wh E_k$ be the event defined exactly as the event $E_R$ in Lemma \ref{lem:exp1}, but for the GFF $h\circ g_m$ instead of $h$, where $g_m(z):= 2^{-m} z$. In other words, $\wh E_k$ is defined exactly as the event $E_R$, except that we consider the disk $\D$ (resp.\ $B_{2^{-k}}(0)$) instead of $B_R(0)$ (resp.\ $\D$). Let $p>0$ and $N\in\N$ be such that Lemma \ref{lem:exp1} holds with $R=2^N$. For any $m\in\N$ let $D_m\subset B_{2^{-m}}(0)$ and $\frk p_m\in\C$ be defined as $D$ and $\frk p$, respectively, in the proof of Proposition \ref{prop:powerlaw}, but for the Gaussian free field $h\circ g_m$ instead of $h$. By applying Lemma \ref{lem:exp1} iteratively, such that we in step $m\in\{0,\dots,k-1\}$ of the iteration condition on $D_{N(k-m)}$ and $\frk p_{N(k-m)}$, we have $\P[\wh E_{Nk}^c]\leq (1-p)^k$. See the proof of Proposition \ref{prop:powerlaw} for a similar argument.
	
	Since $\P[\wh E_{Nk}^c]\leq (1-p)^k$ 
	and $\wh E_k\subset\wh E_{k+1}$ for any $k\in\N$, 
	in order to complete the proof of \eqref{eq:conf1} 
	it is sufficient to show that $\wh E_k\subset E_k$ 
	for any $k\in\N$. If $\wh E_k$  occurs and 
	$z^+\in\BB H\cap (\D\setminus B_{2^{-k}}(0))$ 
	is as in the definition of $\wh E_k$, then the flow line $\eta^L_{z^+}$ (resp.\ $\eta^R_{z^+}$) 
	stays inside $\D\setminus B_{2^{-k}(0)}$ until it hits $\R_-$ (resp.\ $\R_+$). The flow lines $\eta^L_{z_k}$ and $\eta^R_{z_k}$ do not 
	cross the flow lines $\eta^L_{z^+}$ and $\eta^R_{z^+}$, so they stay inside the closure of the domain enclosed by $\eta^L_{z^+}$, $\eta^R_{z^+}$ and $\R$ until they hit $\R$. 
	This implies that $E_k$ occurs, and hence completes the proof of \eqref{eq:conf1}.
	
	Next we will argue that we can find a decreasing sequence $(s_k)_{k\in\N}$ converging to 0 such that
	\eqb
	\P(\wt E_k^c)\leq \exp(-ck/2),\qquad
	\wt E_k:=\big\{ \eta^R_{z_k}([0,t_k^R]) \subset B_{s_k}(0);\,\eta^L_{z_k}([0,t_k^L]) \subset B_{s_k}(0) \big\}.
	\label{eq:conf2}
	\eqe
	By scale invariance of SLE, the probability of $\wt E_k$ is a function of the ratio $|z_k|/s_k$ for fixed $\theta$. Defining $s_k=2^{-\lfloor k/2\rfloor}$, we see by \eqref{eq:conf1} that $\P(\wt E_k^c)=\P( E_{\lceil k/2 \rceil}^c )\leq \exp(-ck/2)$, so \eqref{eq:conf2} holds.
	
	By the Borel-Cantelli lemma, the event $E_k$ occurs for all sufficiently large $k$. By the first characterization of local sets in \cite[Lemma 3.9]{ss-contour} and since flow lines of the GFF are local sets, $x_k\1_{E_k}$ is measurable with respect to $h|_{B_{s_k}(0)}$ for all $k\in\Z$. The sequence $(x_k)_{k\in\N}$ is stationary. So by Birkhoff's Ergodic Theorem,  $\lim_{n\rta\infty} Z_n$ exists a.s.\ and has expectation $p_\theta$. By the Borel-Cantelli lemma, $\lim_{n\rta\infty} Z_n=\lim_{n\rta\infty} \frac{1}{n} \sum_{k=1}^{n}x_k\1_{E_k}$ a.s. Since $\cap_{k\in\N}\sigma(h|_{B_{s_k}(0)})$ is trivial by Lemma \ref{lem-gfftrivial}, we see that this limit is equal to a deterministic constant a.s., so $\lim_{n\rta\infty}Z_n=p_\theta$ a.s.
\end{proof}

\newcommand{\etalchar}[1]{$^{#1}$}

\end{document}